\documentclass[12pt,a4paper,oneside]{article}
\usepackage{amsmath,amstext,amssymb,amsopn,amsthm,color}
\theoremstyle{plain}
\newtheorem{thm}{Theorem}[section]

\newtheorem{lem}[thm]{Lemma}
\newtheorem{prop}[thm]{Proposition}
\newtheorem{cor}[thm]{Corollary}

\theoremstyle{definition}

\theoremstyle{remark}
\newtheorem*{rem*}{Remark}

\newcommand{\R}{\mathbb{R}}
\newcommand{\N}{\mathbb{N}}

\newcommand{\C}{\mathbb{C}}
\renewcommand{\H}{\mathbb{H}}
\newcommand{\ind}{{\bf 1}}
\renewcommand{\leq}{\leqslant}
\renewcommand{\geq}{\geqslant}

\renewcommand{\sp}[2]{\langle #1,#2\rangle}
\newcommand{\pref}[1]{(\ref{#1})}

\def\o{\over}
\def\({\left(}
\def\){\right)}
\def\[{\left[}
\def\]{\right]}
\def\<{\langle}
\def\>{\rangle}

\def\ra{\rightarrow}

\title {Poisson kernels of half--spaces  in real hyperbolic spaces
\footnotetext{2000 MS Classification:
    Primary 60J45, 60J65;
    Secondary 58J65 . {\it Key words and phrases}: hyperbolic spaces, stable
    processes, Poisson kernel.  Research partially supported by KBN
    grant 2 P03A 041 22 and RTN Harmonic Analysis and Related Problems
    contract HPRN-CT-2001-00273-HARP}}
\author{ T. Byczkowski\\ Institute of Mathematics, Wroc\l{}aw
    University of Technology, Poland
  \\and\\ P. Graczyk \\ D\'epartement de
    Math\'ematiques, Universit\'e d'Angers, France
  \\and\\ A. St\'os \\ Institute of Mathematics, Wroc\l{}aw
    University of Technology\\ Laboratoire de Math\'ematique, UMR 6620,
    Universit\'e Blaise Pascal\\ Clermont-Ferrand, France}
\date{}

\begin{document}
\maketitle

\begin{abstract}
  We provide an integral formula for the Poisson kernel of half-spaces for
  Brownian motion in real hyperbolic space $\H^n$. This enables us to
  find asymptotic properties of the kernel.  Our starting point is the
  formula for its Fourier transform.
  When $n=3$, $4$ or $6$ we give an explicit formula for the Poisson
  kernel itself.  In the general case we give various asymptotics and
  show convergence to the  Poisson kernel of $\H^n$.
\end{abstract}

\newpage
\section{Introduction}

 Investigation of the  hiperbolic Brownian motion is an important
 and intensely developed topic in recent years(cf. \cite{revista}, \cite{bougerol}).
  On the other hand, it is well known that the Poisson kernel for a region is a fundamental
tool in  harmonic analysis or probabilistic
potential theory.  In the classical situation of the Laplacian in $\R^n$,
the exact formula for the kernel leads to many
important results concerning behaviour of harmonic functions.
Moreover, probabilistic potential theory uses Poisson kernel
techniques to give solutions to
the Schr\"odinger equation (\cite{ChZ}).
Availability of the exact formula for the kernel is always of crucial
importance for the argument.

The aim of this paper is to give a representation formula for the
Poisson kernel of a half--space in the real hyperbolic space $\H^n$,
i.e.  for the probability distribution of the  hiperbolic Brownian motion stopped
when leaving a half--space,
and to use it in order to prove exact asymptotics of the kernel.
Note that the boundary of the considered half--space is a horocycle
in $\H^n$.

 The Poisson kernel of a half--space is closely related to stable laws and functionals of the
Brownian motion (\cite{BTF}, \cite{BTFY}, \cite{Y1},\cite{Y2}).
Another motivation comes from the  risk
theory in financial mathematics (\cite{D}).
Our kernel, up to a passage from the dimension $2$
 to the dimension  $n$, was identified in terms of its Fourier
transform in \cite{BTF}. It turns out, however, that it is not sufficient for most
(mentioned above) applications.
Unfortunately, a formula for the kernel itself or its
asymptotical behaviour were not
identified (cf. \cite{BTFY}, p. 589).

From the technical point of view,
the main difficulty is that the inverse Fourier transform
(or the Hankel transform) leads to an
integral containing Bessel functions which has oscillatory character,
see \pref{hs11} below.
Moreover, for integrals like \pref{hs11},  Lebesgue's bounded
convergence theorem  is often not applicable (for example when $|y|\to \infty$)
and we are left with a nontrivial problem of obtaining the asymptotics
of the kernel.

 The paper is organized as follows. In Section 2, after
 some preliminaries, we furnish a new proof of the formula
 for the Fourier transform of  the Poisson kernel of a half--space in $\H^n$.
 Our proof, in contrast to the   proof  in \cite{BTF},
 does not use the result of Dufresne (\cite{D}) on the limiting law
 of an exponential functional of the hyperbolic Brownian motion.
 Finally, we give a first integral formula \pref{hs11}
 for the Poisson kernel of a half--space, based on the inverse
 Fourier transform.

 In Section \ref{pkhf}, in Theorem  \ref{rep} we obtain a second integral formula
 for the Poisson kernel of a half--space. This is our main representation formula.
It is
 much more suitable for further applications than \pref{hs11}.  Section 3
 ends with explicit  integral formulas for
the Poisson kernel of a half--space, that arise in
lower dimensions.

 In Section 4 we study
the above mentioned
asymptotics of the Poisson kernel of a half--space in $\H^n$.
We use our main representation formula from Theorem  \ref{rep}
as well as the semigroup and homogeneity  properties of the Poisson kernel.

\section{Preliminaries}

Consider the half-space model of the $n$-dimensional real hyperbolic space
\begin{equation*}
\H^n = \{ (x_1,\ldots,x_{n-1},x_n)\in \R^{n-1}\times \R:\; x_n>0 \}  .
\end{equation*}
The Riemannian metric, the volume element and the Laplace-Beltrami
operator are given by
\begin{equation*}
  ds^2 = {dx_1^2+...+dx_{n-1}^2+dx_n^2 \o x_n^2},
\end{equation*}
\begin{equation*}
  dV={dx_1...dx_{n-1}dx_n \o x_n^n},
\end{equation*}
\begin{equation*}
  \Delta = x_n^2 (\sum_{i=1}^{n} \partial_i^2)- (n-2)x_n\partial_n,
\end{equation*}
respectively (here $\partial_i={\partial \o\partial x_i}$, $i=1,...,n$).

Let $(B_i(t))_{i=1...n}$ be a family of independent classical
Brownian motions on $\R$ with the generator ${d^2\o dx^2}$ (and not
${1\o 2}{d^2 \o dx^2}$) i.e.  the variance $E^0B_i^2(t) = 2t$.
Then the Brownian motion on $\H^n$, $X=(X_i)_{i=1...n}$,
can be described by the following system of stochastic differential
equations
\begin{equation*}
  \left\{
    \begin{array}{ccc}
      dX_1(t) & = & X_n(t)dB_1(t) \\
      dX_2(t) & = & X_n(t)dB_2(t) \\
      . & . &.  \\
      dX_n(t) & = & X_n(t)dB_n(t) - (n-2)X_n(t)dt.
    \end{array} \right.
\end{equation*}
By the It\^o formula one verifies that the generator of the solution
of this system is $\Delta$. Moreover, it can be easily verified that
the solution is given by
\begin{equation*}
  \left\{
    \begin{array}{rcl}
      X_1(t) & = & X_1(0)+\int_0^t X_n(0) \exp(B_n(s)-(n-1)s)dB_1(s) \\
      X_2(t) & = & X_2(0)+\int_0^t X_n(0) \exp(B_n(s)-(n-1)s)dB_2(s)\\
      . & . &.  \\
      X_n(t) & = &X_n(0)\exp( B_n(t)-(n-1)t)
    \end{array} \right.
\end{equation*}

Convention:  by $c$ (or $C$) we always denote a general constant that
depends on $n$ and other constant parameters only. The value of these
constants may change in the same string of estimates.

Below we identify the Poisson kernel
(the function $y \to P_a(x,y)$) in terms of its Fourier transform.
Theorems  \ref{hs0} and \ref{pkf} cover facts that are essentially known.
 Similar results  can be found  in both \cite{BTF}
and \cite{BTFY},  in a slightly different setting of more general generators on
the space $\H^2$.
In order to make  our  paper self-contained,
we  include these facts here, in the present setting of $\H^n$
and  with a different short proof of the formula \pref{hs11}.

Define the projection $\tilde {} : \R^n\ni u=(u_1,...,u_n) \to \tilde
u=(u_1,...,u_{n-1})\in \R^{n-1}$.  In particular,  $\tilde X(t) =
(X_1(t),...,X_{n-1}(t))$.

Consider a half-space $D=\{ x\in \H^n:\; x_n>a \}$ for some fixed
$a>0$. Define
\begin{equation*}
  \tau=\inf\{ t\geq 0:\; X(t)\notin D\} =
\inf\{ t\geq 0:\; X_n(t) = a\}.
\end{equation*}
By $P_a(x,dy)$, $x=(x_1,x_2,...,x_n)\in D$,
$y=(y_1,y_2,...,y_{n-1},a)\in \partial D$ we denote the
Poisson kernel of $D$, ie. the distribution  of $X(\tau)$
starting at $x$ (since $X_n(\tau)=a$ it is enough to consider
the distribution of $\tilde X(\tau)$).

\begin{thm}\label{hs0}
\begin{equation*}
  {\cal F}[P_a(x,\cdot)](u) = E^x \exp(i\< u,\tilde X_\tau\>)
  = \exp(i\sp{\tilde x}{u}) \left( x_n\o a \right)^{n-1 \o2}
  {K_{n-1 \o2}(|u|x_n)\o K_{n-1 \o2}(|u|a)},\quad
    x\in D,\; u\in \R^{n-1}
\end{equation*}
where $K_\nu$, $\nu>0$, is the modified Bessel function of  the third kind,
called also Macdonald function.
\end{thm}
\begin{proof}

Since $B_i(t)$ are independent and $\tau$ depends only on $X_n$
(i.e. $B_n$) we obtain
\begin{equation}
  \label{hs1}
  E^x \exp(i \sp{u}{\tilde X(\tau)})
  =E^{x_n} \prod_{j=1}^{n-1} E^{\tilde x} \exp(i u_jX_j(\tau)).
\end{equation}
We adopt here a useful notation $E^{x_j} Y=E(Y|X_i, i\not=j), E^{\tilde x}Y=
E(Y|X_n)$ for the conditional expectation.
We compute the  integral $\int_0^\tau X_n(s)\;dB_j(s)$ by approximation
 of $\ind_{\{s<\tau \}}(s)X_n(s)$ by  simple processes of the form $\sum_k f_k
 \ind_{[t_k,t_{k+1})}$, where $f_k\in {\cal F}(B_n(t_k))$.
Using independence of the increments of $B_j$ and the fact that the
function under the integral below is independent of $B_j$, we get
\begin{eqnarray*}
  E^{\tilde x} \exp(i u_jX_j(\tau)) &=& E^{x_j} \exp(iu_jX_j(\tau))
  \\&=& E^{x_j}\exp\left(iu_j(X_j(0)+\int_0^\tau X_n\;dB_j)\right)
  \\&=& \lim e^{iu_jx_j}E^{x_j} \exp\left( iu_j \sum_k
    f_k(B_j(t_{k+1})-B_j(t_k))\right)
  \\&=& e^{iu_jx_j}\lim\prod_k
  \exp(-u_j^2f_k^2(t_{k+1}-t_k))
  \\&=&e^{iu_jx_j} \exp\left(
    -u_j^2\int_0^\tau X_n^2(s)\;ds\right).
\end{eqnarray*}
This and (\ref{hs1}) imply that
\begin{equation*}
  E^x \exp(i \sp{u}{\tilde X(\tau)})
  = \exp(i\sp{\tilde x}{u}) E^{x_n}\exp\left(-|u|^2\int_0^\tau
  X_n^2(s)ds\right) =
\exp(i\sp{\tilde x}{u}) E^{x_n}e_q(\tau),
\end{equation*}
where $e_q(\tau) = \exp(\int_0^\tau q(X_n(s)))ds$ with
$q(y) = -(|u|y)^2$.
Observe that the function $\varphi(y) = E^y e_q(\tau)$ is
by definition the {\it gauge} for the Schr\"odinger operator $L+q$
based on the generator $L$ of $X_n$ and the potential $q$.
By general theory (see
e.g. \cite{ChZ}, Prop.4.13, p.119) it is a solution for the Schr\"odinger
equation. Since $dX_n(t) = X_n(t) dB_n(t)
-(n-2) X_n(t)dt$, by a standard argument based on the It\^o formula,
 we get the generator of $X_n$
\begin{equation*}
  L=x_n^2{d^2\o dx_n^2}-(n-2)x_n{d\o dx_n}.
\end{equation*}
Consequently,  $\varphi$ satisfies the following equation
\begin{equation}
  \label{hs2}
  y^2\varphi''(y)-(n-2)y\varphi'(y)-|u|^2y^2\varphi(y) = 0
\end{equation}
on the positive half-line.
Let $\varphi(y)=y^{n-1 \o2} g(y)$. Then
$\varphi'(y) = {n-1 \o 2} y^{n-3\o2} g(y) +y^{n-1 \o2}g'(y)$,
$\varphi''(y) = {(n-1)(n-3) \o 4}y^{n-5\o2}g(y) +(n-1)y^{n-3\o2}g'(y)
+y^{n-1\o2}g''(y) $ and consequently (\ref{hs2}) reads as
\begin{equation*}
  y^2g''(y)+yg'(y)- \left(|u|^2y^2+ \left( (n-1)/ 2\right )^2
  \right)g(y)=0.
\end{equation*}
Substituting $|u|y=z$ and $g(y) = h(z)$ we get
\begin{equation}\label{new1}
  z^2h''(z)+zh'(z)- \left(z^2+ \left( (n-1) /2\right )^2
  \right)h(z)=0.
\end{equation}
This is the modified Bessel equation of order $\nu=(n-1)/2$.
Taking into account the form of the general solution
of (\ref{new1}) we infer that
\begin{equation*}
  \varphi(y) = y^{n-1 \o2} (c_1I_{n-1\o2}(|u|y)+c_2K_{n-1 \o2}(|u|y)),
\end{equation*}
for an appropriate choice of $c_1$ and $c_2$,
where $I_{n-1\o2}(\cdot)$ and $K_{n-1\o2}(\cdot)$ are
the modified Bessel function of the first and third
kind, respectively.
Observe that by definition $\varphi(y)$ is bounded in $y$ and
$\varphi(a)=1$. Since $I_{n-1\o2}(|u|y)$ is unbounded
and  $K_{n-1\o2}(|u|y)$ is bounded when $y\ra\infty$, it follows that $c_1=0$. From the
other condition we get the normalizing constant
\begin{equation*}
  c_2={ 1\o a^{n-1\o2}K_{n-1 \o2}(|u|a) }
\end{equation*}
This completes the proof.
\end{proof}
\begin{rem*}
   As $K_\nu(x)\sim x^{-1/2}e^{-x}$ when $x\ra\infty$,
    the Fourier transform of our kernel is in $L^1$.   Thus  there exists
  the corresponding density which we denote   by $P_a(x,y)$.
\end{rem*}
From Theorem \ref{hs0} it follows that nonzero $\tilde x$  gives rise
just to a translation of $P_a(x,y)$ as a function of $y$. Therefore,
in what follows we may and do
assume $\tilde x =0$. Consequently, we may  simplify the notation by
identifying  $P_a(x,y)=P_a((0,...,0,x),y)$, $x>0$.

For notational convenience, let $s={n\o2}-1$ and $\nu=(n-1)/2$. For
$z>0$ we define  (\cite{GR}, 8.432.8)
\begin{equation}\label{ms}
  m_s(z) = \int_0^\infty e^{-u}u^s(u+2z)^s du = d_n e^z
  z^{n-1\o2}K_{n-1\o2}(z),
\end{equation}
where
$d_n=\pi^{-{1\o 2}} 2^{s+{1\o 2}}\Gamma\left(n\o 2\right)$.
Observe that for $n\in 2\N$ the
function $m_s$ is just a polynomial of the degree $s$. In this case we regard
$m_s(z)$ as defined for all complex numbers.
  By \cite{GR}, (8.468, p. 915), for $n\in 2\N$ we get
  \begin{equation*}
    m_s(z) = d_n \sqrt{\pi\o2}
    \sum_{j=0}^s { (2s-j)! 2^{j-s} \o
      j!(s-j)!}z^j\/.
  \end{equation*}
In particular, $m_0(z)=1$, $m_1(z)=2(1+z)$ and $m_2(z) =8(z^2+3z+3)$.
\begin{thm} [Poisson kernel formula]\label{pkf}
  Let $a>0$, $x>a$ and $y\in \R^{n-1}$. If $|y|>0$ then
  \begin{equation}
    \label{hs11}
    P_a(x,y)= \left( x \o 2\pi a \right)^\nu
   |y|^{-{n-3\o2}} \int_0^\infty   {K_\nu(rx) \o  K_\nu(ra) }
   J_{n-3\o2}(r|y|)r^{n-1\o2}dr,\ \ \ \ \ \nu={n-1\o 2}
  \end{equation}
  and when $|y|=0$ it is understood in the limiting sense, i.e.
  \begin{equation}\label{hs14}
    P_a(x,0) = {2^{2-n}\o\Gamma\(n-1\o2\)}
    \left( x \o \pi a \right)^\nu
   \int_0^\infty   {K_\nu(rx) \o  K_\nu(ra) }
   r^{n-2}dr.
  \end{equation}

  Denoting  $s={n\o2}-1$  we have
  \begin{equation}\label{hs15}
      P_a(x,y) =  (2\pi)^{-\nu} |y|^{-{n-3\o2}} \int_0^\infty
      e^{-r(x-a)}{m_s(rx) \o m_s(ra)} r^{n-1\o2}
      J_{n-3\o2}(r|y|)dr, \quad |y|>0.
  \end{equation}
  The special case $|y|=0$ reads as
  \begin{equation*}
      P_a(x,0) = {2^{2-n}\pi^{-\nu}\o \Gamma\(n-1\o2\)} \int_0^\infty
      e^{-r(x-a)}{m_s(rx) \o m_s(ra)} r^{n-2}dr.
  \end{equation*}
\end{thm}
\begin{proof}
  Recall that if $f$ is a radial function, $f(y)=f_o(|y|)$, then so
  is ${\cal F}f$ and  the Fourier inversion formula   in $\R^{n-1}$
  reads, up to a factor $(2\pi)^{-(n-1)}$,  as the Hankel transform of order $(n-3)/2$ (\cite{F},
  (7.38), p. 247):
  \begin{equation*}
    f_o(|y|) = (2\pi)^{-{n-1\o2} }\int_0^\infty
    ({\cal F}f)_o(r)(r|y|)^{1-{n-1\o2}}J_{{n-1\o2}-1}(r|y|)r^{n-2}dr.
  \end{equation*}
  This gives \pref{hs11}.
  Now, \pref{hs15} is immediate and
  the special cases $y=0$
  follow from the asymptotics  of the Bessel function (see
  e.g. \cite{GR}, 8.440  or \cite{F}, (5.10), p.130)
  \begin{equation}\label{as8}
    J_\nu(z) \sim {1\o2^\nu\Gamma(1+\nu)}z^\nu,\quad z\to 0.
  \end{equation}
  The proof is complete.
\end{proof}
\begin{cor}\label{homogeneity}
The Poisson kernel $P_a(x,y)$, as a function of three variables
$(a,x,y)$, is a homogeneous function of order
$-n+1$:
$$
P_{ta}(tx,ty)=t^{-n+1} P_a(x,y),\ \ \ \ \ t>0.
$$
\end{cor}
\begin{proof}
This is obvious by a change of variables $\tilde r=tr$ in the formula
(\ref{hs11}) written for $P_{ta}(tx,ty)$.
\end{proof}
\begin{rem*}
  Certainly, when $n\in 2\N$ then also $J_{n-3\o2}(r|y|)$ simplifies
  to an elementary function. This fact, however, is not very useful in
  what follows and we will not pursue this further.
\end{rem*}

\section{Poisson kernel of half-space}\label{pkhf}

In this section we give a  representation formula for the Poisson
kernel. For  $n = 2$ the resulting
formula coincides with the one of $\R^2$, so that below we shall
always tacitly assume $n\ge 3$ (note, however, that a great part of our argument
remains valid also for $n=2$).

 From now on we use the following  notation, partially introduced in the preceding Section:
 \begin{eqnarray*}
& \nu=\frac{n-1}{2},\ \ s=\frac{n}{2}-1\\
& \lambda=x-a,\ \  \rho=|y|
 \end{eqnarray*}
and  with a little  abuse
of language we identify $P_a(x,y) = P_a(x,\rho)$.

The following technical lemma is essential in what follows.

\begin{lem}\label{dzielenie}
Let
  \begin{equation*}
    Q(z) = z-(\nu^2-1/4) {\lambda \o 2ax},\quad z\in \C.
  \end{equation*}
 Define  $F_\lambda(z)$ by the following formula:
  \begin{equation} \label{Flambda0}
    \lambda F_\lambda(z) =
    {(z/a) e^{\lambda z/a}(x/a)^{\nu} K_{\nu}(xz/a) - (x/a)^{\nu -{1\o2}}
    Q(z/a) K_{\nu}(z) \o    K_{\nu}(z)}.
  \end{equation}
  Then
  \begin{equation}\label{zanik}
    F_\lambda(z) = O(z^{-1}),\quad  z\to \infty.
  \end{equation}
\end{lem}
\begin{proof}
   Using the asymptotic expansions for the modified Bessel function $K_\nu(z)$
   (\cite{GR}, 8.451.6, p. 910) we get
   \begin{equation}\label{asymp-ms}
    e^z z^\nu K_\nu(z) =
     z^{\nu -{1\o2}} (c_0+{c_1 \o 2z} + R_2),
   \end{equation}
   where
   \begin{equation*}
     c_k=c_k^{(n)}=\sqrt{\pi\o2} {\Gamma({n/2}+k) \o k!
     \Gamma({n/2}-k)},\quad k=0,1,
   \end{equation*}
   $R_2 = O(z^{-2})$ and $|z|$ is large enough.  Hence, it is enough to show
   \begin{equation}\label{reductio}
     (z/a)e^{xz/a}(xz/a)^\nu K_\nu(xz/a)- (x/a)^{\nu-{1\o2}} Q(z/a)e^z
     z^\nu K_\nu(z) =
     O(z^{\nu-{3\o2}}),\quad z\to\infty.
   \end{equation}
   From \pref{asymp-ms} it follows that on one hand we have
   \begin{equation*}
     (z/ a) e^{xz/a} (xz/a)^\nu K_\nu(xz/a) = c_0\(x/ a\)^{\nu-{1\o2}} ({z^{\nu
     +{1\o2}} / a})
     + (c_1 /(2 x))\( x/a\)^{\nu -{1\o2}} z^{\nu-{1\o2}} +O(z^{\nu -{3\o2}}).
   \end{equation*}
   On the other hand, using ${n(n-2)\o 4 }c_0 = c_1$, we get
   \begin{eqnarray*}
     \(x\o a\)^{\nu -{1\o2}} Q\(z\o a\) e^z z^\nu K_\nu(z) & = &
     \(x\o a\)^{\nu -{1\o2}} \( {z\o a} - {n(n-2)\o 8} {\lambda \o ax}\)
     (c_0 z^{\nu-{1\o2}}+{1\o2}c_1z^{\nu -{3\o2}}+ O(z^{\nu-5/2}))
     \\ =
     {c_o \o a} \(x\o a\)^{\nu -{1\o2}} z^{\nu+{1\o2}}
     &+&\(x\o a\)^{\nu -{1\o2}} \({c_1 \o 2a}
     - {n(n-2)\o 4 }c_0 {\lambda\o 2ax} \)z^{\nu -{1\o2}}
     + O(z^{\nu -{3\o2}})
     \\ =
     {c_o \o a} \(x\o a\)^{\nu -{1\o2}} z^{\nu +{1\o2}}
     &+&\(x\o a\)^{\nu -{1\o2}} {c_1 \o 2x} z^{\nu -{1\o2}} + O(z^{\nu -{3\o2}}).
   \end{eqnarray*}
   Then \pref{reductio} is obviously satisfied and the assertion follows.
\end{proof}
\begin{rem*}
  The advantage of this lemma is due to the fact that we may and do
  use it for $z\in\C$. This fact is exploited below.
\end{rem*}
Observe that the function $K_\nu(z)$ has no zeros in $\{\Re
(z) \ge 0\}$ (cf. \cite{E1}, p. 62) and hence $F_\lambda(z)$ is
analytic in this half-plane.
Moreover,
by the inverse Laplace transform theorem (\cite{F}, Theorem 8.5)
together with \pref{zanik}
we get that $F_\lambda$ is the Laplace transform of some
function $w_\lambda$, i.e.
\begin{equation}\label{wlambda}
  F_\lambda(z) = \int_0^\infty e^{-zv}w_\lambda(v)dv\/,
\end{equation}
under the additional condition that for some $b>0$ the following limit
\begin{equation} \label{Laplace}
\lim_{r \to \infty} {1 \o 2\pi i} \int_{b-ir}^{b+ir} F_\lambda(z) e^{vz} dz
\end{equation}
exists for all $v>0$ and it is a piecewise continuous function of $v$
admitting the Laplace transform. Then  the limit is equal
to $w_\lambda(v)$.

The existence of the above limit is shown in Theorem \ref{formula},
together with
an explicit formula for the function $w_\lambda$ itself.

We are ready to state our representation formula.

\begin{thm}\label{rep}
  \begin{equation}\label{rep1}
    P_a(x,\rho) = \frac{\Gamma(s)}{2\pi^{n/2} }{\lambda\o (\lambda^2+\rho^2)^{s+1}}
    \int_0^\infty {w_\lambda(v)L(\lambda,\rho,v) \o
      ((\lambda+av)^2 +  \rho^2)^s} dv,
  \end{equation}
where $L(\lambda,\rho,v)$ is the following function
\begin{equation*}
 L(\lambda,\rho,v) =
 s((\lambda+av)^2-\lambda^2)((\lambda+av)^2+\rho^2)^s
 - (\lambda^2+\rho^2) [((\lambda+av)^2+\rho^2)^s-(\lambda^2+\rho^2)^s].
\end{equation*}
\end{thm}
\begin{proof}
  Recall that (\cite{E}, vol I, (7) and (8) p. 182 or \cite{GR}, 17.13.43,44 )
  \begin{equation*}
    \int_0^\infty e^{-\lambda r} r^\nu J_{\nu-1}(r \rho) dr
    = 2^\nu\pi^{-{1\o2}}\Gamma(\nu + {1\o2})\rho^{\nu-1}
    {\lambda \o (\lambda^2+\rho^2)^{\nu+{1\o2}}}
  \end{equation*}
  and
  \begin{equation*}
    \int_0^\infty e^{-\lambda r} r^{\nu-1} J_{\nu-1}(r \rho) dr
    = 2^{\nu-1}\pi^{-{1\o2}}\Gamma(\nu-{1\o 2})\rho^{\nu-1}
    {1\o (\lambda^2+\rho^2)^{\nu-{1\o2}}}.
  \end{equation*}

For $z=ra$ we have
   \begin{equation}\label{wlambda1}
     \lambda F_\lambda(ra) =
     {r m_s(rx) - \( x/ a\)^s Q(r)m_s(ra) \o m_s(ra) },
   \end{equation}
so
   \begin{equation*}
     {r m_s(rx) \o m_s(ra) }
     = \(x\o a\)^s Q(r) + \lambda F_\lambda(ra)
   \end{equation*}
and hence by \pref{ms}
\begin{equation*}
  \( x\o a \)^{\nu} {r K_\nu(rx) \o K_\nu(ra) }
    = e^{-\lambda r}[\({x\o a}\)^sQ(r)+\lambda F_\lambda(ra)].
\end{equation*}
 Putting this into the Hankel transform formula
 \pref{hs11} and using  \pref{wlambda}  we get
  \begin{equation*}
    2\pi^{n/2}\Gamma(s)^{-1} P_a(x,\rho) = \lambda s  \(x\o a \)^s
    \[ {2\o (\lambda^2 +\rho^2)^{s+1}} -
    {1\o 2xa}{s+1 \o (\lambda^2+\rho^2)^s} \]
    +\lambda \int_0^\infty {w_\lambda(v) \o ((\lambda+av)^2+\rho^2)^s}
    dv.
  \end{equation*}
   Putting $r=0$ in \pref{wlambda1} we get
   \begin{equation} \label{mom0}
     \lambda F_\lambda(0) = \lambda\int_0^\infty w_\lambda(v)dv
     =-\(x/ a\)^s Q(0) = s(s+1)\(x/ a \)^s {\lambda\o 2xa}
   \end{equation}
   so that
   \begin{equation*}
     \lambda(F_\lambda(ra) - F_\lambda(0))m_s(ra)
      = rm_s(rx) - (x/a)^s r m_s(ra).
   \end{equation*}
   Dividing both sides by $r$ and taking the limit $r\to 0$ we obtain
   \begin{equation*}
     \lambda a F_\lambda'(0) = -\lambda a \int_0^\infty v
     w_\lambda(v)dv
     = 1- (x/ a)^s.
   \end{equation*}
   We used the fact  that $vw_\lambda(v)$ allows the Laplace transform
   which is evident from Theorem \ref{formula}.
   Hence
   \begin{equation*}
     \( x/ a \)^s = 1+ \lambda a \int_0^\infty vw_\lambda(v)dv
   \end{equation*}
   Moreover,
   \begin{equation} \label{momgen0}
     \lambda[(F_\lambda(ra) - F_\lambda(0))-raF_\lambda'(0)]m_s(ra)
      = rm_s(rx) - rm_s(ra).
   \end{equation}
   Again, dividing both sides by $(ra)^2$, letting $r\to 0$ and using
   $m_s(0)=m_s'(0)$ we get
   \begin{equation*}
    {\lambda F_\lambda''(0) / 2} = {\lambda/  a^2}
   \end{equation*}
   so that
   \begin{equation} \label{mom1}
     1={a^2 \o 2} \int_0^\infty v^2 w_\lambda(v)dv.
   \end{equation}
   The facts that $ F_\lambda''(0)$ exists and that the function
   $v^2 w_\lambda(v)$ admits the Laplace transform follow
   from Theorem \ref{formula}.
    Consequently we have
   \begin{eqnarray*}
     \(x\o a\)^s &=& 1+ \lambda a \int_0^\infty v w_\lambda(v)dv
     \\&=&
     {a^2\o 2}\int_0^\infty v^2w_\lambda(v)dv
     +\lambda a\int_0^\infty v w_\lambda(v) dv
     \\ &=& {1\o2} \int_0^\infty a v (2\lambda+av)w_\lambda(v)dv.
   \end{eqnarray*}
   Finally, $2\pi^{n/2}\Gamma(s)^{-1}P_a(x,\rho)$ is equal to
   \begin{eqnarray*}
    &{}& \lambda s \int_0^\infty {av(2\lambda+av)w_\lambda(v) \o
     (\lambda^2+\rho^2)^{s+1}}dv
   -\lambda \int_0^\infty {w_\lambda(v) \o (\lambda^2+\rho^2)^s}dv
   +\lambda\int_0^\infty {w_\lambda(v) \o ((\lambda+av)^2+\rho^2)^s} dv
   \\ &=&
    {\lambda\o (\lambda^2+\rho^2)^{s+1}}
       \int_0^\infty {w_\lambda(v)L(\lambda,\rho,v) \o
         ((\lambda+av)^2 +  \rho^2)^s} dv,
   \end{eqnarray*}
   and the assertion follows.
\end{proof}
Below we give a description of the function $w_\lambda$.
The formula depends on the zeros of the function $K_\nu(z)$.
Even if in general the values of these zeros are not
given explicitly,  we are able to prove some important properties
(as boundedeness or asymptotics) of $w_\lambda$, which are essential
 in applications. Moreover, in
lower dimensions we provide explicit formulas as well (see Section 3).

The function $K_\nu(z)$ extends  to an entire function
when $n$ is even and has a holomorphic extension to
$\C\setminus (-\infty,0]$ when $n$ is odd.
Denote the set of
zeros of the function $K_\nu(z)$ by
$Z=\{z_1,...,z_{k_\nu}\} $. We give some needed information about
these zeros (cf. \cite{E1}, p.62).
Recall that in the case of  even dimensions,
the functions
 $m_s(z)$, $s=\nu -1/2$,  are
polynomials of degree $s$.
 They always  have the same zeros as $K_\nu$, so $k_\nu=(n/2)-1$
 when $n\in 2\N$. For $n=2k+1$,
$k_\nu$ is the even number
closest to $(n/2)-1$. In particular, for $n=3$ we have $k_\nu
=k_{1/2}=0$, for $n=5$ and $7$ we have $k_\nu=2$. The functions
$K_\nu$ and $K_{\nu-1}$ have no common zeros.

 In order to describe the function $w_\lambda$ we introduce additional notation.
Let, as before $\nu=(n-1)/2$ and define

\begin{equation} \label{w_1}
w_{1,\/\lambda}(v)= -{(x/a)^\nu \o \lambda a} \sum_{i=1}^{k_\nu}
 {z_i e^{\lambda z_i/a} K_\nu(xz_i/a) \o K_{\nu-1}(z_i)} \/ e^{z_i v}\/,
\end{equation}
and $w_{1}^{\ast}(v) = \sup_{0<\lambda \leq a} |w_{1,\/\lambda}(v)| $

 Using the functions $m_s$,  the formula \pref{w_1} reads as follows:

\begin{equation*} 
w_{1,\/\lambda}(v)={-1 \o (n-2)\lambda a}
\sum_{i=1}^{k_\nu} {m_s(xz_i/a) \o m_{s-1}(z_i)} \/ e^{z_i v}\/.
\end{equation*}

We define additionally in the case of $n$ odd  (so $\nu\in\N$ )
\begin{equation} \label{w_2}
 w_{2,\/\lambda}(v) =(-1)^{\nu+1} {(x/a)^\nu \o \lambda a}
 \int_0^\infty {
   I_\nu\(xu/ a\)K_\nu(u) - I_\nu(u)K_\nu\(xu/ a\)\o
 K_\nu^2(u)+\pi^2I_\nu^2(u)} \/ e^{-\lambda u/a} e^{-vu} \/u du\/,
\end{equation}
and, as before $w_{2}^{\ast}(v)=\sup_{0<\lambda \leq a} |w_{2,\/\lambda}(v)| $.

We also need the following asymptotic formulas for the modified Bessel
functions $I_\nu$,
$K_\nu$: For $u \geq 1$ we have (\cite{F},\cite{GR})
\begin{equation}\label{kinfty}
I_\nu(u)=  (2\pi u)^{-{1\o2}} e^u [1+ E_1(u)] \/, \quad
K_\nu(u)= \pi^{{1\o2}} (2 u)^{-{1\o2}} e^{-u} [1+ E_2(u)] \/ ,
\end{equation}
where $ E_1(u), E_2(u)=O(u^{-1})$, $u\to\infty$.

When $u \to 0$ we have for $\nu>0$:
\begin{equation}\label{kizero}
I_\nu(u) \sim c_\nu u^\nu, \quad
K_\nu(u) \sim c'_\nu u^{-\nu}
\end{equation}
with $c_\nu=2^{-\nu}/\Gamma(\nu+1)$ and
$c'_\nu=2^{\nu-1}\Gamma(\nu)$. Whenever $\nu =0$ one has
$I_0(u) \sim 1$, $K_0(u) \sim \log (2/u)$.
We now formulate and prove our representation theorem for the function
$w_\lambda$.

\begin{thm}\label{formula}
 In the even dimensions
\begin{equation*}
w_\lambda(v) = w_{1,\/ \lambda}(v);
\end{equation*}
while, in the odd dimensions
\begin{equation*} 
w_\lambda(v) = w_{1,\/ \lambda}(v) +  w_{2,\/ \lambda}(v)\/.
\end{equation*}
Moreover, we have $|w_i^{\ast}(v)| \leq C_i(\nu,a)$, $i=1, 2$ and
\begin{equation*}
w_1(v) =  \lim_{x\to a+} w_{1,\/\lambda}(v) = {1 \o a^2}
    \sum_{i=1}^{k_\nu} z_i^2 e^{z_i v}\/;
\end{equation*}
$$
(-1)^{\nu+1}w_{2,\lambda}(v)\geq 0, \ \ \ v\geq 0,\ \ \ \ \  (n\ {\rm odd});
$$
\begin{equation*}
   w_2(v) = \lim_{x\to a+} w_{2,\/\lambda}(v) =
   {(-1)^{\nu+1} \o a^2}
   \int_0^\infty {ue^{-vu} du \o K_\nu^2(u) + \pi^2 I_\nu^2(u)},\/\ \ \ \ \  (n\ {\rm odd});
\end{equation*}\begin{equation*}
   \int_0^\infty v^k w_{1}^{\ast}(v) dv <\infty, ~~ k=1, 2,
   \ldots \/;\quad
   \int_0^\infty v^{n-1} w_{2}^{\ast}(v) dv <\infty \/;
\end{equation*}\begin{equation*}
 \lim_{v \to \infty} v^k w_{1\/,\/\lambda}(v) = 0,\quad  k=1,2,\ldots \/;
\end{equation*}\begin{equation*}
 \lim_{v \to \infty} v^{n+1} w_{2\/,\/\lambda}(v) = {(-1)^{\nu+1}   n!
   \o2^{n-2}\Gamma(\nu)\Gamma(\nu+1)}
   {(x/a)^{n-1}-1 \o \lambda a},\/\ \ \ \ \ (n\ {\rm odd}).
\end{equation*}
\end{thm}
\begin{proof}
  We recall the basic formula \pref{Flambda0}
  \begin{equation*}
   \lambda F_\lambda(z) = { (z/a) e^{\lambda z/a}(x/a)^{\nu}
      K_{\nu}(xz/a)  -  (x/a)^{\nu -{1\o2}} Q(z/a) K_{\nu}(z) \o
      K_{\nu}(z)}
   \end{equation*}

By
 standard rules for
computing residues of meromorphic functions and using the following formula for
derivatives of Bessel functions (cf. \cite{E1}, 7.11(22) p.79)
\begin{equation*}
  {d\o dz} (z^{\nu} K_{\nu}(z)) = -z^\nu K_{\nu-1}(z),
\end{equation*}
we obtain
\begin{equation} \label{res}
\mathrm{Res }_{z_i} F_\lambda =
 -{(x/a)^{\nu} \o \lambda a}  { z_i e^{\lambda z_i /a}
       K_{\nu}(xz_i/a)  \o
      K_{\nu -1}(z_i)}.
 \end{equation}
  Using the functions $m_s$,  we obtain
 \begin{equation} \label{res2k}
 \mathrm{Res }_{z_i} F_\lambda=
 {-1 \o (n-2)\lambda a}
  {m_s(xz_i/a) \o m_{s-1}(z_i)}.
  \end{equation}

As mentioned before (see \pref{wlambda} and \pref{Laplace}), by the
inversion theorem for the Laplace transform we have
  \begin{equation*}
    w_\lambda(v) = {1\o 2\pi i} \lim_{r\to \infty} \int_{b-ir}^{b+ir}
    F_\lambda(z) e^{zv} dz
  \end{equation*}
  for some $b>0$. We show the existence of the above limit together with computing
  formula for the function $w_\lambda$.


  The technique of integration is different in even and odd
  dimensions. This is due to the fact that in the first case the
  function under the integral extends to a meromorphic one while in
  the odd dimension we have to deal with a branch cut.

  For $n\in 2\N$ we choose any $b>0$. All the zeros of
  $m_s(z)$ satisfy $\Re(z_i) < b$ (actually, we have in general
  $\Re(z_i)<0$, $i=1,...,s$, cf. \cite{E1},  p. 62).
  To calculate $w_{\lambda}$ we integrate over the rectangular
  contour with corners at $b-ir$, $b+ir$, $-r-ir$, $-r+ir$. By
  \pref{zanik}
 we infer that integrals over the upper, left, and
  bottom side of the rectangle tend to 0 as $r\to \infty$.
  Hence, by the residue theorem, the limit in \pref{Laplace} exists and
  is equal to the sum of all residues of the
  function $F_{\lambda}(z)e^{z\lambda}$.
  Thus, we have $w_{\lambda}= w_{1,\/\lambda}$ and the assertion follows.

  In the odd dimensions, however, the function under integral is no
  longer meromorphic. We make the branch cut along the negative real
  axis $(-\infty,0]$ and
  change the contour of integration to wrap around this line (see the
  picture).
  
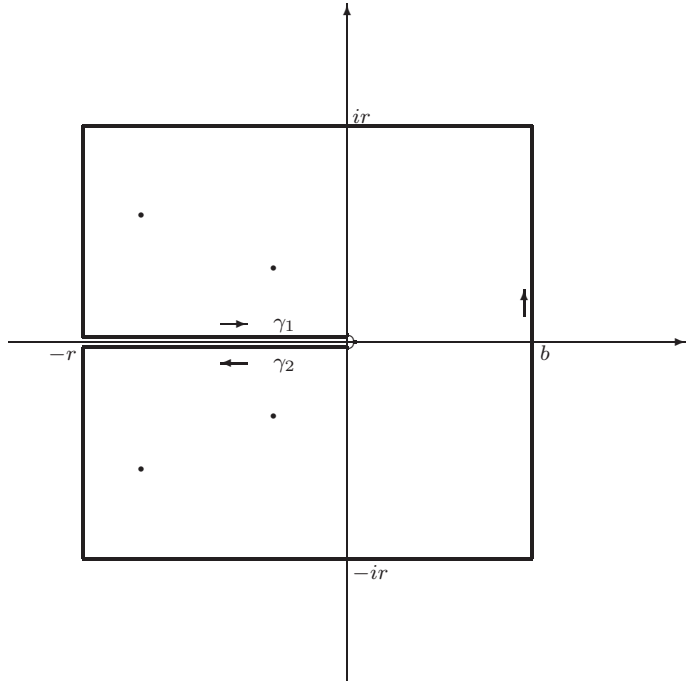
\begin{center}
  \begin{picture}(256,256)(0,0)
    \put(0,128){\vector(1,0){256}}
    \put(128,0){\vector(0,1){256}}
    \linethickness{1pt}
    \put(128,128){\oval(5,5)[r]}
    \put(128,126){\line(-1,0){100}}
    \put(128,130){\line(-1,0){100}}
    \put(28,126){\line(0,-1){80}}
    \put(28,130){\line(0,1){80}}
    \put(28,210){\line(1,0){170}}
    \put(28,46){\line(1,0){170}}
    \put(198,46){\line(0,1){164}}
    \put(201,123){$_b$}\put(130,40){$_{-ir}$}
    \put(130,213){$_{ir}$}\put(15,123){$_{-r}$}
    \linethickness{0.5pt}
    \put(195,138){$\vector(0,1){10}$}
    \put(100,135){$_{\gamma_1}$} \put(100,120){$_{\gamma_2}$}
    \put(80,135){$\vector(1,0){10}$} \put(90,120){$\vector(-1,0){10}$}
    \put(100,100){\circle*{2}}\put(50,80){\circle*{2}}
    \put(100,156){\circle*{2}}\put(50,176){\circle*{2}}
   \end{picture}
   \end{center}
  First, we examine behaviour of our function near the negative axis
  $(-\infty,0)$.  For $z=-y$ ($y>0$) we have (see \cite{E1}, (45), p. 80)
  \begin{equation*}
  \lim_{\epsilon \to 0+}  K_\nu(-y+i\epsilon) = e^{-i \pi \nu}K_\nu(y) -i\pi I_\nu(y),
  \end{equation*}
  \begin{equation*}
    \lim_{\epsilon \to 0+} K_\nu(-y-i\epsilon) = e^{i \pi \nu}K_\nu(y) +i\pi I_\nu(y).
  \end{equation*}

Now, observe that, similarly as before, the integrals over the left,
upper and bottom side of our rectangular contour vanish as
$r\to\infty$ by \pref{zanik}. The same holds true for the
half-circle with radius $\epsilon\to 0$ around the origin.
Note that the branch cut and the residues for $F_\lambda(z)$ are due to
the term
\begin{equation*}
  \tilde F_\lambda(z) = { z(x/a)^\nu e^{\lambda z/a}
      K_{\nu}(xz/a)  \o
      \lambda a K_{\nu}(z)},
\end{equation*}
the rest of the function $F_\lambda(z)$ being holomorphic in $\C$.
Therefore
\begin{equation*}
  {1\o 2\pi i}\int_{b-ir}^{b+ir}  F_\lambda(z)e^{zv} dz
 = w_{1,\/\lambda}(v) -{1\o 2\pi i}\(\int_{\gamma_1}+\int_{\gamma_2}\)\tilde F_\lambda(z)e^{zv} dz.
\end{equation*}

 After taking the limits $r\to \infty$ and   $\epsilon\to 0$, we get
\begin{eqnarray*}
 \(\int_{\gamma_1}+\int_{\gamma_2}\) \tilde F_\lambda(z)e^{zv}dz
 &=&
 {(x/a)^\nu \o \lambda a} \(
- \int_{0}^\infty {u
   (e^{-i\pi\nu}K_\nu(xu/a)-i\pi I_\nu(xu/a)) \o
   e^{-i\pi\nu}K_\nu(u)-i\pi I_\nu(u) } e^{-\lambda u/ a} e^{-vu}
   du \right.
 \\ &+& \left.
 \int_0^\infty {u
   (e^{i\pi\nu}K_\nu(xu/a)+i\pi I_\nu(xu/a)) \o
   e^{i\pi\nu}K_\nu(u)+i\pi I_\nu(u) } e^{-\lambda u/ a} e^{-vu}
   du
\)
\\ &=&(-1)^{\nu+1} {(x/a)^\nu \o \lambda a} 2\pi i
 \int_0^\infty {u
   [I_\nu(u)K_\nu(xu/a)-I_\nu(xu/a)K_\nu(u)]\o
 K_\nu^2(u)+\pi^2I_\nu^2(u)} e^{-\lambda u/ a} e^{-vu}
 du.
\end{eqnarray*}
This ends the proof of the first part of the Theorem.

All that remains is to show the corresponding properties of the functions
$w_{i,\/\lambda}$, $i=1,2$. We begin with $w_{1,\/\lambda}$, which is
easier to analyze.
First of all, observe that $\Re (z_i) <0$, so for fixed $\lambda>0$ the
function $w_{1,\/\lambda}$ is bounded
and $\lim_{v \to \infty} v^k w_{1,\/\lambda}(v)=0$, for all $k=1,2, \ldots $.
To see what happens when $\lambda \to 0$ we use the formula for the residue  of
$F_\lambda$ (see \pref{res}), together with the Lagrange
formula. Since $K_\nu(z_i)=0$ we get
\begin{eqnarray*}
\mathrm{Res }_{z_i} F_\lambda &=&
 {1 \o \lambda}  { (z_i /a) e^{\lambda z_i /a}(xz_i / a)^{\nu}
      K_{\nu}(xz_i/a)  \o
     - z_i^{\nu}K_{\nu -1}(z_i)}
  \\ &=& {- e^{\lambda z_i /a} \o a^2} {z_i \o  z_i^{\nu -1}K_{\nu -1}(z_i)}
      { (xz_i / a)^{\nu} K_{\nu}(xz_i/a)  -z_i^{\nu} K_{\nu}(z_i) \o
      (xz_i/a)-z_i}
   \\ &=& { e^{\lambda z_i /a} \o a^2} {z_i \o  z_i^{\nu -1}K_{\nu -1}(z_i)}
         (\xi z_i)^{\nu} K_{\nu -1}(\xi z_i)
   \\ &\to & {z_i z_i^{\nu}K_{\nu -1}(z_i) \o a^2 z_i^{\nu-1}K_{\nu -1}(z_i)}
   = \({z_i \o a}\)^2\/,
   \end{eqnarray*}
because $1<\xi<x/a$ and  $\xi \to 1$ as $\lambda \to 0$.

Furthermore, for $0<\lambda \leq a$ we have
\begin{equation*}
|\mathrm{Res }_{z_i} F_\lambda| \leq \({|z_i| \o a}\)^2
2^\nu\sup_{1 \leq \xi \leq 2}\left|{  K_{\nu -1}(\xi z_i) \o
 K_{\nu -1}(z_i) }\right| \/.
\end{equation*}

Since $\Re(z_i)<0$, we have obtained that $|w_1^\ast(v)| $is bounded by
a constant $C_1(\nu,a)$ and that
 $w_1^{\ast}$
integrates all powers of $v$.

We now prove the corresponding statements for $w_{2,\/\lambda}$.
Observe that the numerator in \pref{w_2} is equal to
\begin{equation*}
 K_\nu(u)K_\nu(xu/a)\(
{I_\nu(xu/a) \o K_\nu(xu/a) }-
{I_\nu(u) \o K_\nu(u) }\)
\end{equation*}
and hence is positive,  because the function $I_\nu(u)/K_\nu(u)$, $u>0$, is
obviously increasing.

Using the Lagrange formula once again and taking into account
$(z^\nu I_\nu(z))'= z^\nu I_{\nu-1}(z)$(see \cite{E1} 7.11(19) p.79) we obtain for
$0<\lambda \leq a$
\begin{eqnarray*}
  &{}& {1 \o \lambda} {[I_\nu\(xu/ a\)K_\nu(u) - I_\nu(u)K_\nu\(xu/ a\)]\o
    K_\nu^2(u)+\pi^2I_\nu^2(u)}
  \\ &=&
  {\lambda u \o a} {1 \o \lambda}  {[(xu/a)^\nu I_\nu\(xu/ a\)- u^\nu I_\nu(u)]
    \o (xu/a)-u}
  { K_\nu(u) \o
    (xu/a)^\nu  [ K_\nu^2(u)+\pi^2I_\nu^2(u)]}
  \\ &-&
  {\lambda u \o a} {1 \o \lambda}  {[(xu/a)^\nu K_\nu\(xu/ a\)- u^\nu K_\nu(u)]
    \o (xu/a)-u}
  { I_\nu(u) \o
    (xu/a)^\nu  [ K_\nu^2(u)+\pi^2I_\nu^2(u)]}
  \\ & = &
  {u \o a}
  {(\xi_1 u)^\nu I_{\nu-1}\(\xi_1 u\)K_\nu(u) +  (\xi_2 u)^\nu
    I_\nu(u)K_{\nu-1}\(\xi_2 u\) \o
    (xu/a)^\nu [K_\nu^2(u)+\pi^2I_\nu^2(u)]}
  \\ & \to &
  {u \o a}
  {[I_{\nu-1}\(u\)K_\nu(u) + I_\nu(u)K_{\nu-1}\(u\)]\o
    K_\nu^2(u)+\pi^2I_\nu^2(u)}
  = {1 \o a}\/ {1 \o K_\nu^2(u)+\pi^2I_\nu^2(u)}\/,
\end{eqnarray*}
 where for the last equality we used \cite{E1} 7.11(39) p.80 .
 Here the convergence takes place when $\lambda \to 0$ and $0<\xi_1,
 \xi_2 \leq x/a \leq 2$ and $\xi_1, \xi_2 \to 1$ as $\lambda \to
 0$. Thus, we have obtained
 \begin{equation*}
   \lim_{\lambda \to 0}w_{2,\/\lambda}(v)= {(-1)^{\nu+1} \o a^2} \int_0^\infty {u
     e^{-vu} du \o K_{\nu}^2(u)+\pi^2 I_{\nu}^2(u)},
 \end{equation*}
 since the passage to the limit under the integral sign is justified
 by  \pref{major} below.

 Moreover, using the above equations and the asymptotic behavior
 \pref{kinfty} of
 $I_\nu$ and  $K_\nu$, we obtain for $u\geq 1$ and $0<\lambda \le a$
 \begin{eqnarray}
   &{}& {1 \o \lambda} {[I_\nu(xu/a)K_\nu(u) - I_\nu(u)K_\nu(xu/a)]\o
     K_\nu^2(u)+\pi^2I_\nu^2(u)} e^{-xu/a} e^u \nonumber
   \\& \leq & {u \o a}
   {2^\nu [e^{-xu/a} I_{\nu-1}(xu/a) e^u K_\nu(u) + e^{-u} I_\nu(u) e^u
     K_{\nu-1}(u)] \o (x/a)^\nu (K_\nu^2(u)+\pi^2I_\nu^2(u))} \nonumber
   \\ & \leq &
   cu^2{[1+E_1(xu/a)] [1+E_2(u)]+[1+E_1(u)] [1+E_2(u)] \o a \cosh(2u)}
     \leq {Cu^2 \o  \cosh(2u)}. \label{u1}
 \end{eqnarray}
 For $u\leq 1$ we have
 \begin{eqnarray*}
 &{}& {1 \o \lambda} {[I_\nu(xu/a)K_\nu(u) - I_\nu(u)K_\nu(xu/a)]\o
  K_\nu^2(u)+\pi^2I_\nu^2(u)}
  e^{-xu/a} e^u
   \\ & \leq &
    2^\nu{u \o a}
    { I_{\nu-1}(2u) K_\nu(u) +  I_\nu(u) K_{\nu-1}(u) \o
    K_\nu^2(u)+\pi^2I_\nu^2(u)}  e^u
  \end{eqnarray*}
Now, if $\nu-1>0$ (i.e $n>3$), using the asymptotics \pref{kizero}
we obtain that the above expression is bounded from above by
\begin{equation*}
Cu^{2\nu} {1+u^2 \o 1+\pi^2 u^{4\nu}}\le \tilde Cu^{2\nu},\quad u\in(0,1).
\end{equation*}
For $\nu=1$  (i.e. $n=3$) one obtains in fact the same bound
\begin{equation*}
Cu^{2} {1+u^2 \log(2/u) \o 1+u^{4}}\leq C u^{2},\quad u\in(0,1).
\end{equation*}
Thus, we finally  get  for $\nu\geq 1$
\begin{equation*}
w_2^{\ast}(v) = \sup_{0<\lambda\leq a} |w_{2,\/\lambda}(v)| \leq C_1 \int_0^1
u^{2\nu+1}  e^{-vu} du
 + C_2 \int_1^\infty {u^3 e^{-vu} du \o \cosh(2u)}.
\end{equation*}

Now, one easily obtains
\begin{equation} \label{major}
 w_2^{\ast}(v) \leq C(\nu,a) \quad \text{and}\quad
 w_2^{\ast}(v) \leq C_1/v^{2\nu  +2} + C_2 e^{-v},
\end{equation}
and the conclusions concerning the function $ w_2^{\ast}(v) $  follow.

To finish the proof we show the existence and compute the
limit
\begin{equation*} 
\lim_{v\to\infty} v^{2\nu+2} w_{2,\/\lambda}(v)=
\lim_{v\to\infty} v^{n+1} w_{2,\/\lambda}(v).
\end{equation*}
As before, we take into account the expression
under the integral sign in \pref{w_2} multiplied by $v^{2\nu +2} = v^{n+1}$
and, after changing variables $t=vu$ we obtain

\begin{eqnarray}
 &{}&  {v^{2\nu +2}u[I_\nu\(xu/ a\)K_\nu(u) - I_\nu(u)K_\nu\(xu/ a\)]\o
 K_\nu^2(u)+\pi^2I_\nu^2(u)}
 e^{-\lambda u/a} e^{-vu}\/du \nonumber
  \\ &=&  {v^{2\nu}t[I_\nu\(xt/ av\)K_\nu(t/v) - I_\nu(t/v)K_\nu\(xt/av\)]\o
  K_\nu^2(t/v)+\pi^2I_\nu^2(t/v)}
  e^{-\lambda t/av} e^{-t}\/ dt\/.\label{expr1}
\end{eqnarray}
Using the same formulas \pref{kizero} as before,
we obtain that for any fixed $t>0$ the expression above has the
following asymptotics when $v\to\infty$
\begin{eqnarray*}
 & & {v^{2\nu} t c_\nu c_\nu^{'}[(xt/av)^\nu (t/v)^{-\nu} -
 (xt/av)^{-\nu} (t/v)^{\nu}] \o
 (c_\nu^{'})^2 (t/v)^{-2\nu} + \pi^2 c_\nu^2 (t/v)^{2\nu}}
 e^{-\lambda t/av} e^{-t}\/ \/
 \\ &=& {(x/a)^{n-1}-1 \o (x/a)^{\nu}} {c_\nu c_\nu^{'} t^n \o
 (c_\nu^{'})^2 + \pi^2 c_\nu^2 (t/v)^{4\nu}}
 e^{-\lambda t/av} e^{-t}\/ \/
 \\  &\to & {(x/a)^{n-1}-1 \o (x/a)^{\nu}} \({c_\nu \o c_\nu^{'}}\) t^n
  e^{-t}\/ \/ .
\end{eqnarray*}
Moreover, for any fixed $t>0$ and $v$ such that $t<v$ we get that
\pref{expr1} is bounded by $c(x,a,\nu)t^n e^{-t}$.

Now, we write
\begin{equation*}
 v^{2\nu+2} w_{2,\/\lambda}(v) = (-1)^{\nu+1}{(x/a)^\nu \o \lambda a} \(\int_0^v \ldots +
 \int_v^\infty \ldots \).
\end{equation*}
We use (\ref{u1}) with $t/v=u\geq 1$  and  we observe that the expression in the second integral is bounded from above by
\begin{equation*}
t^2(x/a)^{\nu} v^{2\nu-1} {\exp(xt/av-t/v) \o \cosh(2t/v)}
e^{-\lambda t/av} e^{-t} \leq  v^{2\nu-1}(x/a)^{-\nu} t^2e^{-t}\/.
\end{equation*}
Since $v \to \infty$, the second integral tends to $0$, while the first one
converges to the following limit

\begin{equation*}
\lim_{v\to \infty} v^{n+1} w_{2\/,\lambda}(v) = (-1)^{\nu+1}{(x/a)^{\nu} \o \lambda a}
{ (x/a)^{n-1} - 1\o (x/a)^\nu} \({c_\nu\o c_\nu^{'}}\) \int_0^\infty t^n e^{-t} dt
= (-1)^{\nu+1}n! { (x/a)^{n-1}-1 \o \lambda a }{ c_\nu \o c_\nu^{'}}\/.
\end{equation*}
This ends the proof of the theorem.
\end{proof}
{\bf Examples.}
  \ We finish this Section by  writing down explicit integral formulas  for some special
cases. Observe that in $\H^2$ our Poisson kernel is
identical with the Euclidean one. Thus, the simplest nontrivial
situation arises in $\H^3$. Recall that $\lambda=x-a$, $\rho=|y|$.
\begin{cor}
If $n=3$ then
  \begin{equation*}
    w_\lambda(v)={x \o \lambda a^2} \int_0^\infty { I_1(xu/a)
      K_1(u) - K_1(xu/a) I_1(u) \o K_1^2(u)+\pi^2 I_1^2(u)} e^{-u
      \lambda /a} e^{-vu}u\/du
  \end{equation*}
  and
  \begin{equation*}
    L(\lambda,\rho,v)= 2^{-1} a v(2\lambda+av)((\lambda+av)^2
    +\rho^2)^{{1\o2}} -(\lambda^2 + \rho^2)[((\lambda+av)^2
    +\rho^2)^{{1\o2}} - (\lambda^2+\rho^2)^{{1\o2}}],
  \end{equation*}
  hence
  \begin{equation*}
    P_a(x,y)= { \lambda \o 2\pi (\lambda^2+\rho^2)^{{3\o2}}}
    \int_0^\infty {w_\lambda(v) L(\lambda,\rho,v) \o
      ((\lambda+av)^2+\rho^2)^{{1\o2}} }dv.
  \end{equation*}
  If $n=4$ then $w_\lambda(v)=a^{-2} e^{-v}$ and
  \begin{equation*}
    P_a(x,y) = { \lambda \o 2\pi^2(\lambda^2+\rho^2)^2 }
    \int_0^\infty { (2\lambda+av)^2v^2e^{-v} \o (\lambda+av)^2 +\rho^2}dv
  \end{equation*}
  If $n=6$ then
  \begin{equation*}
    w_\lambda(v)= {3 \o a^3} e^{-3v/2}
[(2\lambda
    +a)\cos(\sqrt{3} v/2) +\sqrt{3} a \sin(\sqrt{3} v/2)]
  \end{equation*}
  and
  \begin{equation*}
    L(\lambda,\rho,v)= (a v(2\lambda+av))^2
    [2av(2\lambda+av)+3(\lambda^2+\rho^2)],
  \end{equation*}
  hence
  \begin{equation*}
    P_a(x,y)= { \lambda \o 2\pi^3(\lambda^2+\rho^2)^3} \int_0^\infty
    {w_\lambda(v) L(\lambda,\rho,v) \o ((\lambda+av)^2+\rho^2)^2 }dv
  \end{equation*}
  \end{cor}
\begin{proof}
In the case $n=3$ the function $K_\nu=K_1$ has no zeros
and   we have $w_\lambda(v)=w_{2,\lambda}(v)$.
 We use the Theorems \ref{rep} and \ref{formula}.

  For $n=4$ we have $m_0(x)=1$, $m_1(x)=2(1+x)$. Certainly,
  $w_\lambda(v)=w_{1,\lambda}(v) $ and
  $ L(\lambda,\rho,v)=  (a v(2\lambda+av))^2 $ so all we have to do
  is to find the function  $w_\lambda(v)$. We apply the Theorem \ref{formula} and
  obtain

  \begin{equation*}
 w_\lambda(v) =w_{1,\lambda}(v) = {-1 \o 2\lambda a} {m_1(-x/a) \o m_0(-1)} e^{-v}
  ={(-1/a) 2 (1-x/a) \o 2\lambda} e^{-v} = a^{-2} e^{-v}.
  \end{equation*}

  If $n=6$ then $s=2$ and $m_2(z) =8(z^2+3z+3) =
  8(z-{-3+i\sqrt{3}\o2})(z-{-3-i\sqrt{3}\o2})$.
  Put $z_1=-3/2 +i\sqrt{3}/2$. According to formula \pref{res2k} we obtain

   \begin{eqnarray*}
   \mathrm{Res }_{z_1} F_\lambda &=&
   {-1 \o (6-2)\lambda a}
    {m_2(xz_1/a) \o m_1(z_1)}
  \\ &=& {-1 \o 4\lambda a} {8( (xz_1/a)^2 +3(xz_1/a)+3) \o 2(1+z_1)}
  \\&=&  {-1 \o  \lambda a} ( (xz_1/a)^2 +3(xz_1/a)+3) (1+\overline{z_1})
  \\&=&  {-1 \o  \lambda a} ( (xz_1/a)^2-z_1^2 +3(xz_1/a)-3z_1) (1+\overline{z_1}    )
  \\&=&  {-1 \o   a^2} [ (1+x/a)z_1^2 (1+\overline{z_1}    )   +3z_1(1+\overline{z_1}    )]
  \\ &=& {3 \o 2a^3} [2\lambda + a - i\sqrt{3} a]\/.
    \end{eqnarray*}

  Finally, we have
  \begin{eqnarray*}
 w_\lambda(v) &=& w_{1,\lambda}(v)
 \\ &=& {3e^{-3v/2} \o a^3} \Re\([2\lambda + a - i\sqrt{3} a] e^{i\sqrt{3} v/2}\)
  \\ &=& {3e^{-3v/2} \o a^3} [(2\lambda +a) \cos(\sqrt{3}v/2) + \sqrt{3} a \sin (\sqrt{3}v/2)]\/.
  \end{eqnarray*}
 This completes the case when $n=6$.
\end{proof}

\section{Asymptotic behavior}\label{asymp}
In this section we study the asymptotic behavior of the Poisson
kernel $P_a(x,\rho)$.
The hardest part is to get the asymptotics
for $\rho\to\infty$  (see Theorem \ref{rhoinfty} below).
It is clear that for integrals like \pref{hs11} Lebesgue's bounded
convergence theorem fails.
Another natural approach by a Tauberian theorem (or the Karamata
theory) does not lead to the solution either. On one hand, at the Laplace
transform level, in the required limit we have to
deal with fine cancellations of divergent integrals.
On the other hand, in the basic cases $n=4$ or $n=6$
our representation formula gives almost immediately the required
asymptotics. This leads to the present approach.

Recall that $s=n/2-1$ and $v=(n-1)/2$, $|y|=\rho$.
We shall compare our results to the behaviour of the classical Poisson
kernel of the upper half-space in $\R^n$,
\begin{equation*}
  P_{\R^n}(x,y) = {\Gamma(n/2) \o \pi^{n/2}} { x \o (x^2 + |y|^2)^{n/2} },
\quad x>0, y\in \R^{n-1},
\end{equation*}
and the Poisson kernel  of the {\it
  entire} hyperbolic space in half-space model,
\begin{equation}\label{pglob}
  P_{\H^n}(x,y) =\frac{\Gamma(n-1)}{\pi^{n-1\o 2} \Gamma\left({n-1\o 2}\right)}
   \left( x\o x^2+|y|^2 \right)^{n-1}, \quad
  x>0,\;y\in\R^{n-1},
\end{equation}
see \cite{Guiv}, \cite{Helg}. The constant in the last formula
is easily determined knowing that $\int_{\R^{n-1}}  P_{\H^n}(x,y)dy=1$
on one side and that, on the other side,
$$\int_{\R^{n-1}}   \left( 1\o 1+|y|^2 \right)^{n-1}dy= {2\pi^{n-1\o 2} \o \Gamma\left({n-1\o 2}\right)}
\int_0^\infty {r^{n-2}\o (1+r^2)^{n-1}}dr= {\pi^{n-1\o 2} \o \Gamma\left({n-1\o 2}\right)} B({n-1\o 2},{n-1\o 2})
$$
according to \cite{GR}, 3.194.3.
For example, in the two particular cases of Propositions
  \ref{xinfty} and   \ref{rho-zero} all the three Poisson kernels
  behave in the same way.\\

  The main tools of our study of the asymptotics of $P_a(x,\rho)$
  are the representation formula  \pref{rep1} from Theorem \ref{rep}
   and the semigroup properties of the Poisson kernel. \\

{\bf Semigroup properties of $P_a(x,\cdot)$.}
  By the strong Markov property we obtain the following semigroup
  property of $P_a(x,y)$.
\begin{prop}\label{propertysemigroup}
Let  $b$ be such that $0<a<b<x$ . Then
  \begin{equation} \label{formulasemigroup}
  P_a(x,y) = \int_{\R^{n-1}} P_{b}(x,z) P_a(b,y-z) dz\/,\ \ \ \  y \in \R^{n-1}.
   \end{equation}
  Denoting $P_{a,x}(y)=P_a(x,y)$  we have
  $$ P_{a,x} = P_{a,b}\, *\, P_{b,x}, \ \ \ \ 0<a<b<x,$$
  where $*$ is the usual convolution in  $\R^{n-1}$.
  \end{prop}
  \begin{proof}
  Observe that $\tau_{b} < \tau_a$ so using the strong Markov property we obtain
 for an arbitrary nonnegative and bounded  Borel measurable function $f$ on $\R^{n-1}$:
 \begin{eqnarray*}
 E^x f(X_{\tau_a}) &=& E^x E^{X_{\tau_{b}}}[f(X_{\tau_a})]
  \\ &=&\int_{\R^{n-1}} P_{b}(x,z)\{\int_{\R^{n-1}}  P_a(z,y) f(y) dy \} dz
  \\ &=& \int_{\R^{n-1}} f(y) \{\int_{\R^{n-1}} P_{b}(x,z) P_a(b,y-z) dz\} dy\/.
  \end{eqnarray*}
Thus, we obtain that almost everywhere the following holds
\begin{equation*}
  P_a(x,y) = \int_{\R^{n-1}} P_{b}(x,z) P_a((z,b),y) dz\/,
\end{equation*}
where we denote  $ P_a(x,y)= P_a((0, \ldots ,0,x),y)$
and, according to this notation we also have
$ P_a((z,b),y) =  P_a(b,y-z)$. Since both sides of the above equation
are continuous as a function of $y$, the formula \pref{formulasemigroup} follows.
\end{proof}

{\bf Remark}.The semigroup formula (\ref{formulasemigroup}) holds also for $a=0$,
with $P_{0,x}(y)=P_{\H^n}(x,y)$. This follows from the fact that
as in   \cite{BTF}, $P_{\H^n}$ is the density of $\tilde X_\infty(x)$,
and the proof of the Proposition \ref{propertysemigroup}
still works in this case.

Moreover, when $a\to b$, $a>b$ or when $b\to a$, $b>a$, then
$P_{a,b}\Rightarrow \delta_0$. Consequently, $\{P_{a,b}\}$
is a $2$--parameter continuous probability semigroup.
It means
that $P_{a,b}$ are the densities of the increments ${Y_b - Y_a}$ of a non--homogeneous Levy
process $\{Y_x\}_{0<x<\infty}$, with the distribution of $Y_x$ equal to
 $\tilde X_\infty(x)$   in $\R^{n-1}$.

{\bf Asymptotics when $a\to 0$.}
 When the boundary of the half-space
 in the Euclidean space  $\R^n$ is
 moving away to $-\infty$, the Poisson kernel
 converges to 0. This is not the case in hyperbolic spaces.
 In $\H^n$ we will show the uniform convergence of $P_a(x,\cdot)$
  to the Poisson
 kernel of $\H^n$, given by \pref{pglob}.

 Note that the weak convergence,
 equivalent to the pointwise convergence of Fourier transforms,
 is  simple to see by a probabilistic argument using $X_{\tau_a}\Rightarrow X_\infty$.
 An easy analytic proof of the pointwise convergence of Fourier transforms
   is based on the   Theorem \ref{hs0}, on the asymptotics
 $v^\nu K_\nu(v)\sim 2^{\nu-1}\Gamma(\nu),\ \ v\to 0,$ and on the fact that
 $$
 {\cal F}\left[P_{\H^n}(x,\cdot) \right](u)={1\o 2^{\nu-1}\Gamma(\nu)}
 (x|u|)^\nu K_\nu(x|u|).
 $$
 The last formula follows e.g. from
 \cite{GR} 6.576.7:
\begin{equation}\label{erd}
  \int_0^\infty r^{\mu+\nu+1}J_\mu(r\rho)K_\nu(rx)dr =
  2^{\mu+\nu}\rho^{\mu}x^{\nu}\Gamma(\mu+\nu+1)
  (x^2+\rho^2)^{-\mu-\nu-1},
\end{equation}
with $\mu>\nu-1$, $x>0$.
Putting  $\nu={n-1\o2}$, $\mu={n-3\o2}+\epsilon$ $(\epsilon>0)$,
we have $\mu+\nu+1=n-1+\epsilon$.  Taking limit $\epsilon\to
0$,  by dominated
convergence theorem
we easily extend \pref{erd} to the special case $\mu=\nu-1$.

\begin{prop}\label{global}
 Let $n\in \N$, $n\ge2$. Then for all $y\in\R^{n-1}$
   and  $ x>0$  we have
 \begin{equation*} \lim_{a\to 0} P_a(x,y) = P_{\H^n}(x,y)
\end{equation*}
and the convergence is uniform with respect to $y\in\R^{n-1}$.
\end{prop}
\begin{proof}
We have, by elementary properties  of the convolution,
for any $0<a<b<x$
$$
\|P_{a,x}-P_{0,x}\|_\infty
=\| P_{a,b}*P_{b,x}- P_{0,a}*P_{a,b}*P_{b,x}\|_\infty
\leq
\| P_{0,a}*P_{b,x}-P_{b,x} \|_\infty.
$$
Note that $P_{0,a}*P_{b,x}$ is the action of a probabilistic
operator $T_{0,a}$ with density $P_{0,a}$ on the continuous
function $P_{b,x}$. The function  $P_{b,x}$ is bounded by
Theorem \ref{pkf}.

  The operators $T_{a,x}$ form a continuous 2--parameter semigroup,
so $\lim_{a\to 0}\|T_{0,a} g-g\|_\infty=0$ for any continuous bounded function $g$.
Thus
$$
\| P_{0,a}*P_{b,x}-P_{b,x} \|_\infty\ \to\ 0,\ \ \ a\to 0
$$
and the assertion of the Proposition follows.

A different proof of the Proposition is also possible,
by justifying the passage with $a\to 0$ under the integral
in (\ref{hs11}) and by the Lebesgue bounded convergence
theorem.
\end{proof}

The Proposition \ref{global} implies
the following limit theorem for the hyperbolic
Brownian motion.

\begin{cor}
Let $ {\bf x } =(x_i)_{i=1,\ldots,n}\in \H^n$ and $X_t$
 be the hyperbolic Brownian motion
starting at  ${\bf x}$. Then $X_{\tau_a}$, the
process $X_t$ stopped when first crossing the
hyperplane $\{y_n=a\}$, converges when $a\rightarrow 0$ to a random
variable  $X_\infty$,
concentrated on the  border $\{y_n=0\}$ of $\H^n$
and with the density
\begin{equation*}
P_0({\bf x},y) =
   {\Gamma(n-1) \o \pi^{n-1 \o2}
    \Gamma\left(n-1\o2\right)}
  \left( x_n\o x_n^2+|y-\tilde{\bf x}|^2 \right)^{n-1}
\end{equation*}
where $\tilde{\bf x}=(x_1,\ldots,x_{n-1})$.
 The convergence
of   $X_{\tau_a}$ to   $X_\infty$ is in the sense of uniform
convergence of  the densities of their distributions,
when we project the  hyperplanes $\{y_n=a\}$ on
the  border $\{y_n=0\}$.
\end{cor}

\begin{rem*}
 By Scheffe's theorem, the distributions
of $\tilde  X_{\tau_a}$ converge to the distribution
of   $X_\infty$ in the total variation norm.
\end{rem*}

  {\bf Asymptotics when  $x\to \infty$.}
  The Poisson kernel $P_a(x,y)$ behaves in the same way as the Euclidean Poisson kernel
  and the Poisson kernel of $\H^n$:
\begin{prop}\label{xinfty}
  We have
  \begin{equation*}
    P_a(x,y) \sim cx^{-n+1},\quad x\to \infty.
  \end{equation*}
\end{prop}
\begin{proof}
  First, observe that for $n\ge 2$
  \begin{equation}  \label{as1}
    \(a/ x\)^\nu \le
    e^{r(x-a)}{K_\nu(rx)\o K_{\nu}(ra)}
    \le\(a/ x\)^{{1\o2}}, \quad r>0.
  \end{equation}
  Indeed, since (cf. \cite{GR} 8.432.8)
  \begin{equation*} 
    K_\nu(z)={\Gamma({1\o2}) \o \Gamma(\nu+{1\o2}) }
    \left( {z\o2} \right)^\nu e^{-z}
    \int_0^\infty e^{-zu}(u+2)^{\nu-{1\o2}}u^{\nu-{1\o 2}}du,
  \end{equation*}
  by the change of variable $\bar u=rxu$ ($\bar u=rau$, respectively)
  we get
  \begin{equation}\label{as11}
    e^{r(x-a)}{K_\nu(rx) \o K_\nu(ra)}
    = \( a\o x\) ^{1\o2}
    {\int_0^\infty e^{-u} (u/(rx)+2)^{n-2\o2}u^{n-2\o2}du
      \o \int_0^\infty e^{-u} (u/(ra)+2)^{n-2\o2}u^{n-2\o2}du }.
  \end{equation}
  Since $x>a$ the above quotient of integrals does not exceed 1 and we
  get the upper bound in \pref{as1}.
  Multiplying  the left--hand side of \pref{as11} by
  \begin{equation*}
    1=\(a\o x\)^{n-2\o2}
    {(rx)^{n-2\o2}\o (ra)^{n-2\o2}}
  \end{equation*}
  we get
  \begin{equation}\label{as12}
    e^{r(x-a)}{K_\nu(rx) \o K_\nu(ra)}
    = \( a\o x\) ^{n-1\o2}
    {\int_0^\infty e^{-u} (u+ 2rx)^{n-2\o2}u^{n-2\o2}du
      \o \int_0^\infty e^{-u} (u+2ra)^{n-2\o2}u^{n-2\o2}du }.
  \end{equation}
  Now, $x>a$ implies that the above quotient of the integrals is  greater
  than 1 and the lower bound in \pref{as1} is verified.

 First, we deal with the special case $y=0$.
  By a simple change of variable $rx=t$  in \pref{hs14} and
 by \pref{as12} we get
  \begin{eqnarray*}
    P_a(x,0) &=& c\left(x\o a \right)^{n-1 \o2} {1\o x^{n-1}}
    \int_0^\infty { K_\nu(t) \o K_\nu\(ta/ x\)}
    t^{n-2} dt.
    \\&=& {c\o x^{n-1}}\int_0^\infty
    e^{-t(1-{a/x})}
    {\int_0^\infty e^{-u} (u+ 2t)^{n-2\o2}u^{n-2\o2}du
      \o \int_0^\infty e^{-u} (u+2ta/x)^{n-2\o2}u^{n-2\o2}du }t^{n-2}dt.
  \end{eqnarray*}
  For each $t>0$, when $x$ increases to infinity,  the denominator
  decreases to  $\int_0^\infty e^{-u}
  u^{n-2}du=\Gamma(n-1)$.
  Hence, for $x>2a$ we have
  \begin{equation*}
    {e^{-t(1-{a/ x})}  \o \int_0^\infty e^{-u}
      (u+2ta/x)^{n-2\o2}u^{n-2\o2}du }
    \le {e^{-{t/2}} \o \Gamma(n-1)}.
  \end{equation*}
  Therefore, by bounded convergence
  theorem the assertion for $|y|=0$ follows.

  Now, assume $|y|>0$.
  Recall that
  \begin{equation}\label{as10}
    K_\nu(z) \sim 2^{\nu-1}\Gamma(\nu)z^{-\nu},\quad z\to 0.
  \end{equation}
  By a simple change of variable $rx=z$ in \pref{hs11} we get
  \begin{eqnarray*}
    P_a(x,y) & = & c{|y|^{-{n-3\o2}}\o xa^{n-1 \o 2}}
    \int_0^\infty
    { K_\nu(z) \o K_\nu\left(z{a/ x}\right)}
    J_{n-3\o2}\(z{|y|/ x}\right)  z^{n-1\o2}dz
  \\  &=&
  c x^{-n+1}\int_0^\infty  K_{\nu}(z)
  { x^{\nu} \o (za)^{\nu} K_{\nu}\(za/ x\)}
  { x^{n-3\o2}J_{n-3\o2}\(z|y|/x\) \o (z|y|)^{n-3\o2}}
  z^{3n-5\o2}dz
\end{eqnarray*}
By \pref{as10} and \pref{as8} the two quotients above converge to a
positive constant when $x\to \infty$.
Moreover, the second one remains uniformly bounded in
$z\in(0,\infty)$ and $x>0$.
For $z< {x/ a} $ by \pref{as10} we get
\begin{equation}\label{bc1}
  K_{\nu}(z) z^{3n-5\o2}
  { x^{\nu} \o (za)^{\nu} K_{\nu}\(za\o x\)}
  \le c K_{\nu}(z) z^{3n-5\o2}
\end{equation}
Using \pref{kinfty},
for $z>{x/ a}$ and $x>2a$, say, we get
\begin{eqnarray*}
  K_{\nu}(z) z^{3n-5\o2}
  { x^{\nu} \o (za)^{\nu} K_{\nu}\(za\o x\)}
  &\le &
  c  \left( x\o a \right)^{n\o2} {z^{3n-5\o2}\o z^{n-1\o2}}
  \exp(-z(1-{a/ x}))
  \\& \le &
  cz^{3n -4\o2}\exp(-{z/2}).
\end{eqnarray*}
By this and \pref{bc1} bounded convergence theorem applies.
Consequently, the whole integral above tends to a positive constant as
$x\to \infty$. The assertion follows.
\end{proof}

{\bf Asymptotics when  $x\to a$.}
The asymptotics below are easy to obtain.
\begin{prop}\label{rho-zero}
  \begin{equation*}
    P_a(x,0) \sim {2^{2-n}\Gamma(n-1)\pi^{-{n-1\o2}} \o
      \Gamma\(n-1\o2\) } (x-a)^{-n+1}, \quad x\to a^+.
  \end{equation*}
\end{prop}
\begin{proof}
  From \pref{as1} it follows that
  \begin{eqnarray*}
    \Gamma(n-1)\(a/ x\)^{\nu}(x-a)^{-n+1}
    &=&    \(a/ x\)^{\nu}
    \int_0^\infty e^{-r(x-a)}r^{n-2}dr
    \\&\le&
    \int_0^\infty
    {K_{\nu}(rx) \o K_{\nu}(ra)}r^{n-2}dr
    \\&\le&
    \(a\o x\)^{{1\o2}}
    \int_0^\infty e^{-r(x-a)}r^{n-2}dr
    \\&=&
    \Gamma(n-1)\(a/ x\)^{1\o2}(x-a)^{-n+1}.
  \end{eqnarray*}
  Combining this and \pref{hs14} completes the proof.
\end{proof}
Much more is required, however, to obtain the following {\it
  Euclidean-like} asymptotics
\begin{equation*}
  P_a(x,\rho) \sim c(x-a),\quad x\to a^+,\ \ \ \rho\not=0.
  \end{equation*}
The justification of this important result is postponed after the proof of
Theorem \ref{rhoinfty}.

{\bf Asymptotics when $\rho\to \infty$.}
The most important and difficult thing to prove is what
happens when $\rho\to \infty$.  By $n$
we denote, as before, the dimension of the considered hyperbolic space
$\H^n$. We assume throughout this section that $n>2$. Recall that $s=n/2 -1$.
Let us rewrite the basic formula for $ P_a(x,\rho)$,  using some notation
more suitable for calculations. Denote
\begin{eqnarray*}
&z=\lambda^2+\rho^2,\\
&\kappa=(\lambda+av)^2-\lambda^2,\\
&u=\kappa/z,\\
&\Phi(u)=(1+u)^{-s}-1+su,\ \ \ \ \ u\geq 0.
\end{eqnarray*}
We then have $(\lambda+av)^2+\rho^2=\kappa+z$, so
\begin{equation*}  
 L(\lambda,\rho,v)=L(\kappa,z)=s\kappa(\kappa+z)^s-z[(\kappa+z)^s-z^s]
  \end{equation*}
 and
  \begin{equation*}
{ L(\kappa,z) \o z(\kappa+z)^s}={s\kappa\o z} -
 1+(\frac{z}{\kappa+z})^s=(1+u)^{-s}-1+su=\Phi(u).
   \end{equation*}
  Consequently,
 \begin{equation} \label{formPa0}
  P_a(x,\rho)={\Gamma(s)\o 2\pi^{n/2}} {\lambda\o (\lambda^2+\rho^2)^s} \int_0^\infty
  w_\lambda(v) \Phi(u)\,dv,
   \end{equation}
   where $u=u(x,a,\rho,v)$.
   Writing the last formula  in the form
    \begin{equation} \label{formPa}
    P_a(x,\rho)={\Gamma(s)\o 2\pi^{n/2}} \rho^{-n+2} {\lambda\o( (\lambda/\rho)^2+ 1)^s} \int_0^\infty
  w_\lambda(v) \Phi(u)\,dv,
    \end{equation}
    we see that in order to get the asymptotics of $P_a(x,\rho)$ when $\rho\to\infty$,
    it is sufficient to obtain the asymptotics of $\int_0^\infty
  w_\lambda(v) \Phi(u)\,dv$ when $\rho\to\infty$.

In the sequel we use  the following standard properties of the oscillating binomial
series
$$
\sum_{j=0}^\infty  (-1)^j {(s)_j\o j!} u^j,\ \ \ \ \ (s)_j=s(s+1)\ldots(s+j-1),
 $$
 related to the function $\Phi$. For all $u\geq 0$ and $l\geq 2$ we have
 \begin{equation}\label{Phi1}
 \left| \Phi(u)-\sum_{2\leq j\leq l-1} (-1)^j{(s)_j\o j!} u^j\right| \leq {(s)_l\o l!} u^l.
 \end{equation}
 Moreover, when $u\to 0+$,
 \begin{equation}\label{Phi2}
 \lim_{u\to 0+} u^{-l} \left[\Phi(u)-\sum_{2\leq j\leq l-1} (-1)^j{(s)_j\o j!} u^j\right]=(-1)^l {(s)_l\o l!}
\end{equation}
We also have  for all $u\geq 0$
 \begin{equation}\label{Phi3}
(-1)^{l} \left[\Phi(u)-\sum_{2\leq j\leq l-1} (-1)^j{(s)_j\o j!} u^j\right] \geq 0.
\end{equation}

Note that for $l=2$ the formulas above involve the function $\Phi$ alone (summation is
performed over empty set of indices).

\begin{lem}\label{katymoment}
Let $2\leq l\leq [n/2]$. Then
$$
\lim_{\rho\to\infty} \rho^{2l} \int_0^\infty w_\lambda(v)
 \left[\Phi(u)-\sum_{2\leq j\leq l-1} (-1)^j{(s)_j\o j!} u^j\right]dv =
 (-1)^l{(s)_l\o l!}\int_0^\infty \kappa^l w_\lambda(v) dv.
 $$
 When $n$ is even, this is true for any $l\in\N$, $l\geq 2$.
\end{lem}
\begin{proof}
 We multiply and divide by $u^l$
the expression  under the integral in the last formula.
When $\rho\to\infty$ then $u\to 0$, so the formula (\ref{Phi2})
applies to
$$
u^{-l} \left[\Phi(u)-\sum_{2\leq j\leq l-1} (-1)^j{(s)_j\o j!} u^j\right].
$$
On the other hand,  $\rho^{2l} u^l \to \kappa^l$. The passage to the limit
under the integral sign is justified by (\ref{Phi1})
and the fact that
\begin{eqnarray*}
&&\int_0^\infty \kappa^{[n/2]} |w_\lambda(v)| dv<\infty,\ \ \ n\  {\rm odd},\\
&&\int_0^\infty \kappa^{m} |w_\lambda(v)| dv<\infty,\ \ \ m\in\N,\  n\   {\rm even}.\\
\end{eqnarray*}
Actually, according to Theorem \ref{formula},
 when $n$ is even, the function $w_\lambda$
has finite moments of all orders and when $n=2k+1$, $k=[n/2]$,
the function $v^{n-1}w_\lambda(v)\sim \kappa^k  w_\lambda(v)$
is integrable.
\end{proof}

The observation contained in the next lemma is crucial for our purposes.
 \begin{lem}\label{powers}
  For hyperbolic spaces of  even dimension $n=2k>4$
  we have for all $\lambda >0$
  \begin{equation}   \label{c^2}
   \int_0^\infty \kappa^j  w_\lambda(v) dv = 0
  \end{equation}
  when $j=2,\ldots,k-1$.

For hyperbolic spaces of  odd dimension $n=2k+1>3$
  we have for all $\lambda >0$
  \begin{equation}   \label{c^3}
    \int_0^\infty \kappa^j  w_\lambda(v) dv = 0
   \end{equation}
 when $j=2,\ldots, k$.
  \end{lem}
\begin{proof}
Consider first the case $n=2k>4$ and suppose that the assertion is false.
Let $j_o\leq k-1$ be the smallest power such that
(\ref{c^2}) does not hold.  By Lemma \ref{katymoment} and the formula
(\ref{formPa}), it follows that there exists
$$
\lim_{\rho\to\infty} \rho^{n-2+2j_o} P_a(x,\rho)>0.
$$
We will show that this is contradictory   with the existence of a
 finite $\lim_{\rho\to\infty} \rho^{2n-2}P_0(x,\rho)$. By the semigroup property
 proved in Proposition \ref{propertysemigroup} (see Remark below its proof) we have
 $$
 P_0(x,\cdot)=P_0(a,\cdot)*P_a(x,\cdot).
 $$
 It follows that for $|y|>M>0$
 \begin{eqnarray*}
 |y|^{2n-2} P_0(x,y)&\geq& \left( {|y|\o |y|+1}\right)^{2n-2}
 \int_{|z|\leq 1} P_0(a,z)P_a(x,y-z)|y-z|^{2n-2} dz\\
&=&\left( {|y|\o |y|+1}\right)^{2n-2} \int_{|z|\leq 1} P_0(a,z)P_a(x,y-z)|y-z|^{n-2+2j_o}|y-z|^{n-2j_o}  dz\\
&\geq& c\left( {|y|\o |y|+1}\right)^{2n-2}   \int_{|z|\leq 1} P_0(a,z)|y-z|^{n-2j_o}  dz\\
&\geq& c_1 (|y|-1)|^{n-2j_o}\to \infty,\ \ \ |y|\to\infty,
 \end{eqnarray*}
 because $n-2j_o\geq 2$.

 In the case $n=2k+1>3$, let us  remark that Lemma \ref{katymoment}
 applies for $l\leq k$. We proceed exactly in the same way as in the
 proof in the   case of $n$ even, with the only difference that now $j_o\leq k$,
 so $n-2j_o\geq 1$ and the final contradiction with $\lim_{\rho\to\infty} \rho^{2n-2}P_0(x,\rho)<\infty$
also  holds.
 \end{proof}

  We now state and prove the main result of this section.

 \begin{thm} \label{rhoinfty}
 We have
   \begin{equation*}
     P_a(x,\rho) \sim c\rho^{-2n+2 },\quad \rho\to \infty.
   \end{equation*}
 \end{thm}
 \begin{proof}
 {\it Case $n=2k$}.
  We   will show  that the formula (\ref{c^2}) fails for
 $j=k$, i.e. that
  \begin{equation}   \label{not0}
   \int_0^\infty \kappa^k  w_\lambda(v) dv \not = 0.
  \end{equation}
  The formula (\ref{not0}) together with Lemmas \ref{katymoment} and \ref{powers}
  imply that $\int_0^\infty w_\lambda(v)\Phi(u)dv\sim c\rho^{-2k}=c\rho^{-n}$, $c\not=0$,
   when $\rho\to\infty$.
  Taking into account  the formula (\ref{formPa}) we obtain the desired result.
  Observe that when $n=4$, the function $w_\lambda(v)=a^{-2} e^{-v}$ is positive so
  the formula (\ref{not0}) is apparent. In this case we only need Lemma \ref{katymoment}
  (Lemma \ref{powers} is not available for this case).  \\
  {\it  Proof of the formula (\ref{not0}).}  Suppose that  (\ref{not0}) is not true,
  i.e.
  \begin{equation*}
   \int_0^\infty \kappa^k  w_\lambda(v) dv = 0.
  \end{equation*}
 Then, using  Lemma \ref{katymoment}
  for $l=k+1$,
  it follows that there exists the limit
  $$
  \lim_{\rho\to\infty} \rho^{2k+2}\int_0^\infty w_\lambda(v)\Phi(u) dv=L
 $$
 and, by (\ref{formPa}),
\begin{equation}\label{ro2n}
 \lim_{\rho\to\infty} \rho^{2n} P_a(x,\rho)={\Gamma(s)\o 2\pi^{n/2}}L <\infty.
  \end{equation}
  We will show that the existence of the limit (\ref{ro2n}) leads to a contradiction
  with the convergence of the Poisson kernels $P_a(x,y)$ to $P_0(x,y)$, $a\to 0+$,
  established in Proposition \ref{global}.

  By the homogeneity property of the kernel $P_a(x,y)$, proved in Corollary \ref{homogeneity},
  we have
  $$
  P_{a^2/x}(a,a\rho/x)=(x/a)^{-n+1} P_a(x,\rho),
  $$
  so we   also have
 \begin{equation}\label{ro2nbis}
  \lim_{\rho\to\infty} \rho^{2n} P_{a^2/x}(a,a\rho/x)<\infty.
 \end{equation}
 We will prove that
 \begin{equation}\label{ro2n3}
 \overline{\lim_{\rho\to\infty}}\rho^{2n} P_{a^2/x}(x,\rho)<\infty.
 \end{equation}
  Set $a_1=a^2/x$. As $0<a_1<a<x$, the semigroup property implies
 that
 $$
 |y|^{2n} P_{a_1}(x,y)=\int_{\R^{n-1}}  P_{a_1}(a,z) |y|^{2n} P_a(x,z-y) dz.
 $$
 We divide the last integral into $\int_{2|z|\geq |y|}+\int_{2|z|<|y|}$ and estimate separately
 both integrals.

 By (\ref{ro2nbis}), we obtain that $\lim_{\rho\to\infty}\rho^{2n}  P_{a_1} (a,\rho)<\infty$.
 This is used in the estimate
 \begin{eqnarray*}
 \int_{2|z|\geq |y|} P_{a_1}(a,z) |y|^{2n} P_a(x,z-y) dz
&\leq &2^{2n} \int_{2|z|\geq |y|} P_{a_1}(a,z) |z|^{2n} P_a(x,z-y) dz\\
 &\leq& c 2^{2n} \int_{2|z|\geq |y|} P_a(x,z-y) dz \leq c 2^{2n}
 \end{eqnarray*}
 where the constant $c$ is common for all $|y|>M>0$. Next, observe
 that if $|y|>2|z|$ then $|y-z|>|y|/2$, so that, using (\ref{ro2n}),
  \begin{eqnarray*}
  \int_{2|z|< |y|} P_{a_1}(a,z) |y|^{2n} P_a(x,z-y) dz
  &\leq& 2^{2n}  \int_{2|z|< |y|} P_{a_1}(a,z) |z-y|^{2n} P_a(x,z-y) dz\\
  &\leq &c 2^{2n}   \int_{2|z|< |y|} P_{a_1}(a,z) dz\leq c 2^{2n}
  \end{eqnarray*}
  and (\ref{ro2n3}) is proved. Note that in  order to prove it, a weaker  hypothesis
  \begin{equation*}
 \overline{\lim_{\rho\to\infty}} \rho^{2n} P_a(x,\rho)<\infty.
  \end{equation*}
  is sufficient. Consequently, denoting $q=a/x$ and iterating the last argument
  $j$ times we get
   \begin{equation}\label{ro2nj}
 \overline{\lim_{\rho\to\infty}} \rho^{2n} P_{q^ja}(x,\rho)<\infty,\ \ \ \ j\in\N.
  \end{equation}
  We denote $a_j=q^ja$, $\lambda_j=x-a_j$, $\kappa_j=(\lambda_j+a_jv)^2-\lambda_j^2$
  and we denote $w_{\lambda_j}(v)$ the function appearing in the representation formula
  \pref{rep1} for the kernel $P_{a_j}(x,\rho)$. The formula \pref{ro2nj} implies that
  for all $j\in\N$
   \begin{equation}\label{0j}
   \int_0^\infty \kappa_j^{n/2} w_{\lambda_j}(v) dv=0
 \end{equation}
 (otherwise $\rho^{2n-2}P_{a_j}(x,\rho)$ converges to a  positive constant when $\rho\to\infty$,
 so   $\rho^{2n}P_{a_j}(x,\rho)$ diverges to $+\infty$).\\
 Recall that
 $$
 w_{\lambda_j}(v)=
w_{1,\/\lambda_j}(v)= \sum_{i=1}^{k_\nu}( {\rm Res}_{z_i} F_{\lambda_j} )e^{z_i v}\/=
 -{(x/a_j)^\nu \o \lambda a} \sum_{i=1}^{k_\nu}
 {z_i e^{x z_i/a_j} K_\nu(xz_i/a_j) \o e^{z_i} K_{\nu-1}(z_i)} \/ e^{z_i v}\/.
 $$
 Writing
 $$
 {\rm Res}_{z_i} F_{\lambda_j} =-{1\o a^{n/2}} {z_i^{1/2} \o \sqrt 2\lambda_j }{x^{n/2}\o e^{z_i} K_{\nu-1}(z_i)}
 \left(2{xz_i\o a_j}\right)^{1/2} e^{xz_i/a_j} K_\nu\left({xz_i\o a_j}\right)
 $$
 and using the asymptotics $ K_\nu(u)= \pi^{{1\o2}} (2 u)^{-{1\o2}} e^{-u} [1+ E(u)]$, $u\geq 1$,
 $E(u)=O(u^{-1})$ when $u\to\infty$, we get for $a\to 0+$
 $$
 a^{n/2} {\rm Res}_{z_i} F_{\lambda_j} =-\left({\pi\o 2}\right)^{1/2} {z_i^{1/2}\o x-a} {x^{n/2}\o  e^{z_i} K_{\nu-1}(z_i)}
 [1+E({xz_i\o a_j})]
 $$
 and
 $$
 \tilde w(v):=\lim_{a\to 0} a^{n/2} w_\lambda(v)=-\left({\pi\o 2}\right)^{1/2} x^{n/2-1}
  \sum_{i=1}^{k_\nu} {z_i^{1/2}e^{z_i v}\o e^{z_i} K_{\nu-1}(z_i)}.
  $$
  Thus we get
  \begin{equation}\label{wtilde}
  \lim_j \kappa_j^{n/2} w_{\lambda_j}(v)=(2x)^{n/2} v^{n/2} \tilde w(v).
  \end{equation}
  By the Lebesgue dominated convergence theorem and by \pref{0j}
  $$
  \lim_j \int_0^\infty \kappa_j^{n/2} w_{\lambda_j}(v) dv= (2x)^{n/2}\int_0^\infty  v^{n/2} \tilde w(v) dv=0.
  $$
  Now, when $j\to\infty$ and $\rho$ is fixed,
   one has $\kappa_j\to 0$ and $u_j={\kappa_j\o \lambda_j^2 +\rho^2}\to 0$. By \pref{Phi2}
   and by  \pref{wtilde}, we get
   \begin{eqnarray*}
   \lim_j\int_0^\infty w_{\lambda_j}(v)\Phi(u_j) dv&=& \lim_j\int_0^\infty w_{\lambda_j}(v)[\Phi(u_j) -
   \sum_{2\leq j\leq k-1} (-1)^j {(s)_j\o j!} u^j_j] dv\\
   &=&(-1)^{n/2}  {(s)_{n/2}\o (n/2)!} \left( {2x\o x^2+\rho^2}\right)^{n/2}
    \int_0^\infty v^{n/2} \tilde w(v) dv=0.
   \end{eqnarray*}
   This implies that $\lim_j P_{a_j}(x,y)=0$ which is false because
   $\lim_j P_{a_j}(x,y)=P_0(x,y)\not=0$.
 The proof in the case $n=2k$ is completed.\\
  {\it Case $n=2k+1$}.
 
  By Lemma \ref{powers}, for $n>3$
  \begin{eqnarray*}
 && \int_0^\infty w_\lambda (v)\Phi(u)dv=
   \int_0^\infty w_\lambda (v)[\Phi(u)-\sum_{2\leq j\leq k} (-1)^j{(s)_j\o j!} u^j]dv.
  \end{eqnarray*}
  This formula is also trivially  true for $n=3$ (there is no sum
on the right--hand side).
  We want to show that $ \int_0^\infty w_\lambda (v)\Phi(u)dv\sim C\rho^{-n}$ when $\rho\to\infty$,
  for a positive constant $C$.
  We divide the last integral into the sum of three integrals
 \begin{equation*}
 \int_0^\infty  =
 \int_0^{A}  +
 \int_A^{\varepsilon z^{{1\o2}}}   +
\int_{\varepsilon z^{{1\o2}}}^\infty ,
     \end{equation*}
where $A>0$ is so big that $(-1)^{\nu+1} w_{2,\lambda}(v)> {c\o v^{n+1}}$
for a  constant $c>0$; this is possible because, by Theorem \ref{formula},
  $(-1)^{\nu+1} v^{n+1} w_{2,\lambda}(v) $
converges when $v\to\infty$ to a positive constant. The value of $\epsilon$
will be chosen in the sequel. Recall that $\nu={n-1\o 2}=k$.

We estimate the integral $\int_0^A$ using the bound \pref{Phi1} for
$l=k+1$ and $|w_\lambda(v)|\leq  c/v^{n+1}$. We obtain
$$
 \int_0^{A}  \left|  w_\lambda (v)[\Phi(u)-\sum_{2\leq j\leq k} (-1)^j{(s)_j\o j!} u^j]\right| dv
 \leq {c\o z^{k+1}}\int_0^{A} {\kappa^{k+1} \o v^{n+1}} dv
 \leq {c_1\o \rho^{2(k+1)}} \int_0^A
 {v^{2(k+1)} \o v^{n+1}} dv = {c_1 A\o \rho^{n+1}}.
 $$
 Next, observe that, again by \pref{Phi1} for
$l=k+1$,
 \begin{eqnarray*}
 \int_A^\infty\left|w_{1,\lambda}(v)[\Phi(u)-\sum_{2\leq j\leq k} (-1)^j{(s)_j\o j!} u^j]\right|  dv
 &\leq&{c\o z^{k+1}}\int_A^\infty |w_{1,\lambda}(v)| \kappa^{k+1} dv\\
 &\leq &{c\o \rho^{n+1}} \int_A^\infty |w_{1,\lambda}(v)| v^{n+1} dv
 \leq  {c_2\o \rho^{n+1}},
\end{eqnarray*}
 since the function $|w_{1,\lambda}(v)|$ decreases exponentially when $v\to\infty$.
 Thus, in order to show that  $ \int_0^\infty w_\lambda (v)\Phi(u)dv\sim C\rho^{-n}$ when $\rho\to\infty$,
 it suffices to prove that
 $$
 \int_A^\infty w_{2,\lambda}(v)[\Phi(u)-\sum_{2 \leq j\leq k} (-1)^j{(s)_j\o j!} u^j] dv\sim C\rho^{-n}.
 $$
 Note that by \pref{Phi3} with $l=k+1$ and by the fact that $(-1)^{\nu+1}w_{2,\lambda}\geq 0$,
 proved in Theorem \ref{formula}, the integrand
 $$ w_{2,\lambda}(v)[\Phi(u)-\sum_{2\leq j\leq k} (-1)^j{(s)_j\o j!} u^j]
 $$
 is non--negative. Moreover, by \pref{Phi2}, for $\epsilon>0$ sufficiently small
 and $v<\epsilon z^{1/2}$,  so that $v^2/z<\epsilon^2$,  $v/z\leq v/\sqrt z<\epsilon$
 and $u<b(x,a)\epsilon$, we have
 $$|\Phi(u)-\sum_{2\leq j\leq k} (-1)^j{(s)_j\o j!} u^j| > c u^{k+1},$$
 for a positive constant $c$.
 On the other hand, by \pref{Phi1}, we have for all $u\geq 0$
  $$|\Phi(u)-\sum_{2\leq j\leq k} (-1)^j{(s)_j\o j!} u^j| < c' u^{k+1}.$$
  We consider $\rho$ so big that $A<\epsilon  z^{1/2}$.
  The last  two estimates imply
  \begin{eqnarray*}
  {c\o z^{k+1}} \int_A^{\epsilon z^{1/2}} {\kappa^{k+1}\o v^{n+1}}dv \leq
  \int_A^{\epsilon z^{1/2}}w_{2,\lambda}(v)[\Phi(u)-\sum_{2\leq j\leq k} (-1)^j{(s)_j\o j!} u^j] dv
  \leq {c'\o z^{k+1}} \int_0^{\epsilon z^{1/2}} {\kappa^{k+1}\o v^{n+1}}dv,
  \end{eqnarray*}
  which, using $\kappa=av(2\lambda+av)$, implies
  \begin{eqnarray*}
  {c_1\o z^{k+1}}\int_A^{\epsilon z^{1/2}} dv\leq
  \int_A^{\epsilon z^{1/2}}w_{2,\lambda}(v)[\Phi(u)-\sum_{2 \leq j\leq k} (-1)^j{(s)_j\o j!} u^j] dv
  \leq {c_2\o z^{k+1} }\int_0^{\epsilon z^{1/2}} dv
  \end{eqnarray*}
  and
  \begin{eqnarray*}
  {c_1\epsilon\o z^{k+(1/2)}}-{c_1 A\o z^{k+1}} \leq
  \int_A^{\epsilon z^{1/2}}w_{2,\lambda}(v)[\Phi(u)-\sum_{2\leq j\leq k} (-1)^j{(s)_j\o j!} u^j] dv
  \leq {c_2\epsilon \o z^{k+(1/2)} }
  \end{eqnarray*}
  \begin{eqnarray*}
  {c_3\o \rho^n}- {c_4\o  \rho^{n+1} }
 \leq \int_A^{\epsilon z^{1/2}}w_{2,\lambda}(v)[\Phi(u)-\sum_{2\leq j\leq k} (-1)^j{(s)_j\o j!} u^j] dv
 \leq {c_5  \o \rho^n  }.
  \end{eqnarray*}
  Thus
 \begin{eqnarray*}
  \int_A^{\epsilon z^{1/2}}w_{2,\lambda}(v)[\Phi(u)-\sum_{2\leq j\leq k} (-1)^j{(s)_j\o j!} u^j] dv\sim
 C \rho^{-n},
 \end{eqnarray*}
 for $C>0$, $\rho\to\infty$.

 In the last integral $\int_{\epsilon z^{1/2}}^\infty w_{2,\lambda}(v)[\Phi(u)-\sum_{2\leq j\leq k} (-1)^j{(s)_j\o j!} u^j] dv$,
 the variable
 $u$ is separated from 0, so the expression $|\Phi(u)-\sum_{2\leq j\leq k} (-1)^j{(s)_j\o j!} u^j|$
 is estimated from above, up to a positive factor, by the highest level term $u^k$.  Consequently
\begin{eqnarray*}
 0< \int_{\epsilon z^{1/2}}^\infty w_{2,\lambda}(v)[\Phi(u)-\sum_{2\leq j\leq k} (-1)^j{(s)_j\o j!} u^j] dv
 &\leq& {c\o z^k} \int_{\epsilon z^{1/2}}^\infty  {\kappa^k\o v^{n+1}} dv\\
 &\leq&  {c'\o z^k} \int_{\epsilon z^{1/2}}^\infty  { v^{n-1}\o v^{n+1} }dv= {c"\o z^k} z^{-1/2}\sim
 {C_1\o \rho^n}.
  \end{eqnarray*}
  Since both integrals
  \begin{eqnarray*}
  \int_A^{\epsilon z^{1/2}}w_{2,\lambda}(v)[\Phi(u)-\sum_{2\leq j\leq k} (-1)^j{(s)_j\o j!} u^j] dv&\sim&
 C \rho^{-n}\\
  \int_{\epsilon z^{1/2}}^\infty w_{2,\lambda}(v)[\Phi(u)-\sum_{2\leq j\leq k} (-1)^j{(s)_j\o j!} u^j] dv&\leq&
    C_1\rho^{-n}
  \end{eqnarray*}
  are positive, we conclude that
 $$
 \int_A^\infty w_{2,\lambda}(v)[\Phi(u)-\sum_{2\leq j\leq k} (-1)^j{(s)_j\o j!} u^j] dv\sim C\rho^{-n}
 $$
 and the proof in the case $n=2k+1$ is finished.
 \end{proof}

We now deal with the remaining asymptotics of $P_a(x,\rho)$
 near the boundary:  $x\to a^+$, $\rho\not=0.$

\begin{thm}\label{xtoa}
  We have
  \begin{equation*}
    P_a(x,\rho) \sim c(x-a),\quad x\to a^+,\ \ \ \rho\not=0.
  \end{equation*}
\end{thm}
\begin{proof}

 We recall the basic formula \pref{formPa0} for $P_a(x,\rho)$

 \begin{equation*} 
   P_a(x,\rho)={\Gamma(s)\o 2\pi^{n/2}} {\lambda\o (\lambda^2+\rho^2)^s} \int_0^\infty
   w_\lambda(v) \Phi(u)\,dv,
    \end{equation*}

We apply here
notation and terminology introduced before the formulation of Theorem \ref{rhoinfty}.
In particular, we have
for $0<\lambda \leq a$
\begin{equation*}
0 \leq \Phi(u) \leq
s {\kappa \o z} =s {av(2\lambda + av) \o {\lambda^2 + \rho^2}} \leq \({a \o v}\)^2 v(2+v)\/.
 \end{equation*}
We also have
\begin{equation*}
\lim_{x \to a+} \Phi(u)=s \({av \o \rho}\)^2 - 1 +
\({\rho^2 \o \rho^2 + (av)^2}\)^s\/,
\end{equation*}
and $\lim_{x \to a+} w_\lambda(v) = w(v) = w_1(v) + w_2(v)\/. $

Using properties of $w_{i,\/\lambda}$ stated in Theorem \ref{formula}
one obtains
\begin{equation*}
\lim_{x \to a+}{P_a(x,\rho)\o\lambda}=
{ (4\pi)^s\Gamma(s) \o \rho^{2s}}
\int_0^\infty w(v)[
s \({av \o \rho}\)^2 - 1 +
\({\rho^2 \o \rho^2 + (av)^2}\)^s]\/dv\/.
\end{equation*}
Using formulas \pref{mom0} and \pref{mom1} from the proof of  Theorem \ref{rep}
one computes
\begin{equation*}
1={a^2 \o 2} \int_0^\infty v^2 w_\lambda(v)\/ dv \to {a^2 \o 2}   \int_0^\infty v^2 w(v)\/ dv
\end{equation*}
and
\begin{equation*}
 \int_0^\infty w(v)\/ dv =
\lim_{x \to a+}  \int_0^\infty w_\lambda(v)\/ dv =
\lim_{x \to a+} s(s+1)(x/a)^s {1 \o 2xa} =
 {s(s+1) \o 2a^2}\/.
\end{equation*}
We thus obtain
\begin{equation*}
\lim_{x \to a+}{ P_a(x,\rho)\o\lambda}=
C[{2s \o \rho^{2s+2}} - {s(s+1) \o 2a^2 \rho^{2s}} +
\int_0^\infty {w(v)\/dv \o (\rho^2 + (av)^2)^s}\/].
\end{equation*}
Observe now that when we multiply the right-hand side of the above equation
by $\rho^{2s}$ and let $\rho \to 0$ then the first term tends to infinity while
the second one is constant and the third one converges to $0$:
\begin{equation*}
\rho^{2s} \int_0^\infty {w(v)\/dv \o (\rho^2 + (av)^2)^s} dv=
 \int_0^\infty {w(v)\/dv \o (1 + (av/\rho)^2)^s} dv \to 0\/.
 \end{equation*}
This observation completes the proof.
\end{proof}

\section*{Acknowledgements}
The authors would like to thank K. Bogdan, T. Kulczycki,
M. Ryznar and T. Zak for stimulating conversations on the
subject. Moreover, the first and the third named authors acknowledge
the hospitality of Universit\'e d'Angers where a substantial part of
the work was done. The third named author wishes to express his gratitude to the MAPMO
laboratory at Orl\'eans for giving the opportunity to work further on this
subject during his post--doctoral fellowship.

\end{document}